\documentclass[11pt,a4paper]{article}

\usepackage{amsmath}
\usepackage{amsfonts}
\usepackage{amssymb}
\usepackage{mathtools}
\usepackage{graphicx}
\usepackage{color}
\usepackage{float}
\usepackage{fullpage}
\usepackage{dsfont}
\usepackage{epsfig}
\usepackage{setspace}
\usepackage{multicolumn}
\usepackage{longtable}
\usepackage{caption}
\usepackage{algcompatible}
\usepackage{booktabs}
\usepackage{url}
\usepackage[natbib=true, style=authoryear, backend=biber, abbreviate=false, mincitenames=2, maxcitenames=3, maxbibnames=999, giveninits=true, uniquename=init]{biblatex}
\renewbibmacro{in:}{}
\renewbibmacro*{volume+number+eid}{%
\printfield{volume}%
\printfield{number}%
\setunit{\bibeidpunct}%
\printfield{eid}}
\DeclareFieldFormat[article]{number}{\mkbibparens{#1}}
\newcommand{\ignore}[1]{}
\def\defi{\vcentcolon=}

\newcommand{\new}{\textcolor{blue}}

\begin{document}

\date{}

\title{Policies for the Operation of an Ambulance Fleet under Uncertainty based on a New Preparedness Metric}

\maketitle
\vspace*{-2cm}

\begin{center}
\begin{tabular}{ccc}
\begin{tabular}{c}
Vincent Guigues\\
School of Applied Mathematics, FGV\\
Praia de Botafogo, Rio de Janeiro  22250-900,  Brazil\\
{\tt vincent.guigues@fgv.br}
\end{tabular}&
&
\begin{tabular}{c}
Anton J. Kleywegt\\
Georgia Institute of Technology\\
Atlanta, Georgia 30332-0205, USA\\
{\tt anton@isye.gatech.edu}\\
\end{tabular}
\end{tabular}
\end{center}

\begin{center}
\begin{tabular}{c}
Victor Hugo Nascimento*\\
School of Applied Mathematics, FGV\\
Praia de Botafogo, Rio de Janeiro  22250-900, Brazil\\
{\tt victorhugo.vhrn@gmail.com}\\
\end{tabular}
\end{center}

\begin{abstract}
Two important decisions in the management of an ambulance fleet are ambulance selection decisions and ambulance reassignment decisions.
Ambulance selection decisions determine what to do when an emergency call arrives (such as choosing what ambulance to dispatch to the emergency or putting the emergency in a queue of emergencies waiting for an ambulance to be dispatched).
Ambulance reassignment decisions determine where to send an ambulance next when it has finished service for an emergency.
Making good ambulance selection decisions and ambulance reassignment decisions is challenging because a decision made at a point in time affects the ability of the emergency medical service to respond to future emergencies (that are typically not known when the decision is made).
We propose a new preparedness metric that quantifies the ability of the emergency medical service to respond to future emergencies.
The preparedness metric can be used to make ambulance selection decisions and ambulance reassignment decisions by solving a tractable optimization problem each time that a decision has to be made.
We compare the performance of the resulting method with $9$~methods that have been proposed in the literature, based on data from a real emergency medical service for a large city.
\end{abstract}

\par {\textbf{Practitioner Summary.}} Emergency medical services must promptly decide which ambulance to dispatch when a call arrives and where to reposition ambulances after completing a service. These decisions are difficult because choosing an ambulance for the current call can reduce the system’s ability to respond quickly to future emergencies.

This study introduces a new preparedness metric that helps quantify how ready the ambulance fleet is to handle upcoming, yet unknown, emergency calls. By incorporating this metric into a fast optimization-based decision rule, dispatchers can make better real-time decisions about both ambulance assignment and post-service repositioning.

Using data from a large urban emergency medical service, we show that the proposed approach leads to more balanced fleet availability and improved system readiness compared to existing dispatch and reassignment policies. The method can be implemented as a decision-support tool to improve operational performance.\\


\section{Introduction}
\label{sec:intro}

This paper proposes a method for making ambulance dispatch decisions that takes into account the consequences of current dispatch decisions on the ability of the Emergency Medical Service (EMS) to respond to future emergencies.
To understand such consequences, it is important to have a basic knowledge of EMS operations, and therefore we give an overview of EMS operations in Section~\ref{sec:EMS operations}.
Then, in Section~\ref{sec:EMS decisions}, we describe the ambulance dispatch decisions considered in this paper.
Section~\ref{sec:literature review} reviews the literature that is closely related to the work in this paper, after which we describe the contributions of this paper.

\subsection{Overview of Emergency Medical Service Operations}
\label{sec:EMS operations}

In this section we give an overview of EMS operations.
It is important to keep in mind that there is great variability in the way things are done in different EMS systems, and thus this description can give only a broad overview of EMS operations, while the practitioner has to take into account many additional operational details.
For more details about EMS dispatch operations, see for example \citet{claw:15,coon:15}.

An EMS operates a variety of ambulances with various crew members.
Different ambulance/crew combinations have different capabilities, and it is important to take these capabilities into account when choosing which ambulance to dispatch to a particular emergency.
An EMS also operates ambulance stations where crew members and ambulances can wait to be dispatched.
Stations can vary from parking spots with minimal additional features, to well equipped facilities where crew members can eat, relax, and sleep, and where ambulances can undergo thorough cleaning and maintenance.
An EMS also has a call center that receives emergency calls.
The organization that operates the call center and the organization that operates the ambulances may be the same or different, and in the USA either can be part of a local government or can be outsourced to a private contractor.
The Emergency Medical Dispatchers (EMDs) at a call center answer the calls and ask the callers a sequence of scripted questions.
More specifically, the questions form a tree that branches according to the answers received.
Two widely used systems of questions are the Medical Priority Dispatch System (MPDS), and the Association of Public-Safety Communications Officials (APCO) system.
These systems classify emergencies based on the caller's answers into 30--40 chief complaint types, each of which is further subdivided into 4--17 subtypes that determine the combination of ambulance/crew capability needed and response time urgency.
After an EMD has classified an emergency, the EMD dispatches ambulances if deemed appropriate (or instructs someone else to dispatch ambulances), and gives pre-arrival instructions to the caller.
Any ambulance can be dispatched to an emergency, even if the ambulance is on its way to another emergency or on its way to an ambulance station, but typically an ambulance that is busy with an emergency is not dispatched to another emergency.
Not all calls require an ambulance to be dispatched (calls that do not require an ambulance to be dispatched are not included in our model).
Also, an emergency may be placed in queue, that is, an ambulance is not dispatched immediately (within a minute or two) after the call has been received, but rather, the dispatcher waits until a sufficient number of ambulances are available before an ambulance is dispatched to the emergency location.
Thus, placement of an emergency in queue may be needed if no ambulances are available when the emergency call is received, or may be desirable if the emergency is not urgent and there are only a few available ambulances, so that it is preferred to have the ambulances available for possibly more urgent future emergencies.

When an ambulance arrives at an emergency location, the crew members perform the tasks at hand according to their training.
If one or more patients need transport to a hospital, then the ambulance(s) and crew transport these patients to the hospital deemed best under the circumstances.
The best hospital depends on both the location of the emergency as well as the type of emergency.
Sometimes permission is obtained from the hospital's emergency department before the patients are taken to the hospital, and a hospital sometimes denies permission (``diverts'' the ambulance), but permission is not always asked.
After an ambulance arrives with a patient at a hospital, the patient is transferred to the hospital's emergency department.
The transfer may be quick, or may take a substantial amount of time, depending on decisions made by the emergency department.
After the patient has been transferred, the ambulance crew have to clean the ambulance and complete their report.
If the ambulance requires relatively light cleaning, then the crew can do the cleaning on the spot, and thereafter the ambulance is ready to be dispatched to an emergency in queue, or to be sent to an ambulance station.
If the ambulance requires major cleaning, then the ambulance is taken to an ambulance station where such major cleaning can be done, and thereafter the ambulance is ready to be dispatched to an emergency in queue, or it waits at the ambulance station to be dispatched.

\subsection{Decisions for Ambulance Fleet Operations}
\label{sec:EMS decisions}

The management of ambulance operations requires decisions to be made quickly, with important consequences for the mortality and morbidity of people.
More specifically, as part of the management of a fleet of ambulances, the following decisions must be made:
\begin{enumerate}
\item
When a call requesting emergency aid arrives, a decision has to be made whether to dispatch an ambulance to the emergency location immediately, or whether to place the request in a queue.
If it is decided to dispatch an ambulance to the emergency location, then it also has to be decided which ambulance to dispatch.
We call this decision the ambulance selection decision.
Sometimes it is called a call-initiated decision.
\item
When an ambulance finishes its task with an emergency (for example, at the emergency location or at a hospital after dropping off patients), then it has to be decided what the ambulance should do next.
If there are emergencies waiting in queue, and it is decided to send the ambulance to an emergency waiting in queue, then it also has to be decided to which emergency to send the ambulance.
If the ambulance is not sent to an emergency in queue, then it has to be decided to which ambulance station to send the ambulance, where it can wait for its next dispatch.
We call this decision the ambulance reassignment decision.
Sometimes it is called an ambulance-initiated decision.
\end{enumerate}

In its response to an emergency, an ambulance goes through some or all of the following steps as part of its service:
\begin{itemize}
\item[(a)] travels to the emergency scene;
\item[(b)] provides service at the emergency scene;
\item[(c)] transports the patient(s) to a hospital;
\item[(d)] stays at hospital waiting for the patient to be transferred;
\item[(e)] travels to a cleaning station to clean the ambulance;
\item[(f)] cleans the ambulance at the cleaning station;
\item[(g)] travels to an ambulance station for staging.
\end{itemize}
An ambulance attending an emergency may transport the patient(s) to a hospital or not, and may need to go to a cleaning station after service or not.

Ambulance management decisions involve important trade-offs.
For example, an ambulance may be dispatched to a current, less urgent, emergency, and in the process it may not be available for a future, more urgent, emergency.
These trade-offs are challenging for various reasons.
First, current emergencies are known (at least the location and something about the nature of the emergency is usually known), whereas typically future emergencies are not known.
However, it is known that certain types of emergencies tend to occur with greater frequencies in specific parts of the city and during specific times of the week, and therefore it is prudent to send more available ambulances to these parts during these times.
For example, penetrating trauma and traffic incidents tend to occur with greater frequency in some parts of the city on Friday evenings.
Second, different ambulances and crews have different capabilities to improve the outcomes for patients.
For example, many Emergency Medical Services have basic life support (BLS) and advanced life support (ALS) ambulances.
Some also have other types of ambulances, such as intermediate life support ambulances, stroke units, motorcycles, and helicopters.
Different crew members also have different capabilities and qualifications, including Emergency Medical Technicians (EMTs), Advanced Emergency Medical Technicians (AEMTs), paramedics, and physicians with various specialties.
Third, the consequences of response time and ambulance/crew capabilities are different for different types of emergencies, and for many types of emergencies the relationship between response time or ambulance/crew capabilities and patient outcome is not yet well known.
For example, it is well known that for cardiac arrest the response time is crucial and is more important than the advanced capabilities of the ambulance and crew.
For some emergencies it is known that advanced capabilities, such as ability to administer intravenous treatment or specific pharmaceuticals, are more important.
And for many emergencies there is a trade-off between response time and ambulance/crew capabilities, so that given a choice of dispatching a BLS ambulance that is 10 minutes from the emergency or an ALS ambulance that is 20 minutes from the emergency, either may be better than the other depending on the type of emergency.
Our models make provision for these trade-offs.

\subsection{Literature Review}
\label{sec:literature review}

Many ambulance related optimization problems have been considered in the literature.
Most of these problems address the location of ambulance stations or the assignment of ambulances to stations.
Surveys of this literature can be found in \citet{swer:94,gree:04,galv:08,ingo:13,arin:17,reut:17,bela:19}.
Less research has been done on ambulance selection and reassignment decisions.
This paper focuses on ambulance selection and reassignment decisions, and therefore we describe only the related literature.

\subsubsection{Different approaches to taking consequences of ambulance dispatch decisions on future operations into account}

Among the different approaches to taking the consequences of ambulance dispatch decisions on future operations into account, two extremes are (1)~ignoring the consequences of ambulance dispatch decisions on future operations, and (2)~using detailed models of emergency requests and ambulance operations to forecast the consequences.
The first approach is simple and tractable, but sometimes results in decisions with adverse consequences, whereas the second approach requires great computational effort and usually results in better decisions.
An intermediate approach is to devise a tractable approximation of the consequences, and to test whether it provides better decisions than the first approach.

\paragraph{Approaches that do not model consequences.}

For the ambulance selection decision, the most popular policy in the literature that ignores the consequences is the closest available ambulance rule.
As the name indicates, this policy dispatches the available ambulance that is closest to the emergency.
This policy has been used in the work of \citet{hend:99,hend:04,maxw:09,maxw:10,maxw:13,alan:13}.
A number of other static ambulance selection policies have been proposed.
For example, \citet{band:14} proposed the following policy for an EMS with two priority levels of emergencies.
Ambulances are assigned to home stations with one ambulance per station, only ambulances at their home stations are available for dispatch, and after completing service an ambulance returns to its home station.
For high priority emergencies, the closest available ambulance is selected.
Ambulances are ordered according to the fraction of emergencies that they would serve if all emergencies were served according to the closest ambulance rule, and for low priority emergencies, the least busy available ambulance is selected.
\citet{mayo:13} compared the following related policies.
The EMS service region is partitioned into zones, and contiguous zones are combined into districts such that each district contains at least one ambulance station.
Policy~1 used the closest available ambulance rule within each district, and ambulances cannot serve emergencies in other districts; policy~2 used the closest available ambulance rule within each district, and if no ambulance is available in a district, then the heuristic of \citet{band:14} was used to select an ambulance from another district; policy~3 used the the heuristic of \citet{band:14} to select an ambulance within each district, and ambulances cannot serve emergencies in other districts; and policy~4 used the the heuristic of \citet{band:14} to select an ambulance within each district, and if no ambulance is available in a district, then the heuristic of \citet{band:14} was used to select an ambulance from another district.
Similarly, \citet{lisay:16} compared the closest available ambulance rule with a policy that dispatches the closest available ambulance to high priority emergencies, and the ambulance within a specified response time radius which has the least utilization to low priority emergencies.

Recall that the ambulance reassignment decision may assign an ambulance that has just become available to an emergency in queue, or may send the ambulance to a location to wait for its next dispatch.
\citet{band:12,mayo:13,band:14,lisay:16,jagt:17b} ignored the queueing of emergency requests (it is assumed that emergencies that arrive when no ambulances are available are ``transferred to another service'').
\citet{ande:07} made provision for queueing of emergency requests if no ambulances are available, but did not specify how the decision was made which emergency in queue to serve next when an ambulance becomes available.
\citet{lees:11} dispatched the ambulance to the emergency in queue with location that is closest to the ambulance (called the nearest-neighbor policy), whereas \citet{schm:12,jagt:17a} dispatched the ambulance to the oldest emergency in queue (called the first-come-first-served policy).
To choose a location for an ambulance when it becomes available and is not dispatched to an emergency in queue, a simple policy is to keep the ambulance in place (for example, at the hospital where the ambulance delivered a patient).
\citet{ande:07,lees:11} considered such ambulance reassignment policies.
A popular ambulance reassignment policy in the literature assigns in advance a home station to each ambulance, and when an ambulance becomes available and is not dispatched to an emergency in queue, then the ambulance returns to its home station.
Such ambulance reassignment policies were considered in \citet{gold:90b,hend:04,rest:09,band:12,knig:12b,maso:13,mayo:13,band:14,lisay:16,jagt:17a,jagt:17b}.

\paragraph{Approaches that model and optimize consequences.}

\citet{schm:12} considered a dynamic programming model of EMS operations.
The state space of the dynamic program was too large to compute an optimal policy, and therefore an approximate dynamic programming method was used to choose ambulance selection decisions and to choose a location for an ambulance when it becomes available and is not dispatched to an emergency in queue.
It was shown that the resulting policy consistently outperforms the policy that uses the closest available ambulance rule for ambulance selection decisions, and that returns an ambulance to its home station when it becomes available and is not dispatched to an emergency in queue.
\citet{band:12} also formulated an ambulance dispatching problem as a continuous-time Markov Decision Process (MDP).
The MDP can be solved if the number of zones and number of ambulances are sufficiently small, in which case the resulting policy outperforms the closest available ambulance rule.
\citet{nasr:18} proposed a number of basis functions to approximate the value function of the MDP.
The response time performance under the resulting policy was compared with a number of benchmark policies and with the lower bounds of \citet{maxw:14}.
\citet{guiklevhn2022} proposed a detailed multistage model of EMS operations.
Instead of solving the intractable multistage stochastic integer program, each time an ambulance selection decision or ambulance reassignment decision has to be made, a two-stage stochastic program is solved with the current decision in the first stage and sample paths of continuous relaxations of future decisions in the second stage.
In numerical tests the resulting policy outperformed $6$~alternative policies from the literature.

\paragraph{Approaches that approximate consequences.}

Approaches that ignore operational consequences often perform poorly, and approaches based on detailed models of operations often take too much time to provide recommendations to decision makers.
Therefore methods have been proposed that approximate consequences, while providing recommendations in less time than approaches based on detailed models of operations.
This is the approach pursued in this paper.
\citet{ande:07} described such a method based on a metric of ``preparedness'' that measures how well available ambulances can respond to expected future emergencies.
Their preparedness metric is computed as follows.
For each zone~$\ell \in \mathcal{L}$, identify the $A_{\ell}$ (currently available) ambulances that contribute most to preparedness in zone~$\ell$.
Index these ambulances $1,\ldots,A_{\ell}$ such that the travel times $t^{a}_{\ell}$ of ambulances indexed~$a$ from their current locations to zone~$\ell$ satisfy $t^{1}_{\ell} \le t^{2}_{\ell} \le \cdots \le t^{A_{\ell}}_{\ell}$.
Let parameter $\gamma^{a}$ denote the contribution of ambulance~$a$ to the preparedness in zone~$\ell$.
Ambulances with less travel times are more likely to be dispatched, and therefore $\gamma^{1} > \gamma^{2} > \cdots > \gamma^{A_{\ell}}$.
Let $\lambda_{\ell}$ be a weight proportional to the demand in zone~$\ell$.
Then the preparedness metric $\psi_{\ell}$ of zone~$\ell$ is given by
\[
\psi_{\ell} \ \ \defi \ \ \frac{1}{\lambda_{\ell}} \sum_{a=1}^{A_{\ell}} \frac{\gamma^{a}}{t^{a}_{\ell}}
\]
that is, by the ambulance ``availability'' $\sum_{a=1}^{A_{\ell}} \gamma^{a} / t^{a}_{\ell}$ relative to the demand $\lambda_{\ell}$.
For any zone~$\ell \in \mathcal{L}$ and available ambulance $a$, let $\psi^{-a}_{\ell}$ denote the preparedness metric of zone~$\ell$ after removing ambulance~$a$.
The proposed preparedness metric was used in the following policy.
For ambulance selection decisions, for the highest priority emergencies, the closest available ambulance rule is used.
For lower priority emergencies, the ambulance with expected travel time less than a specified threshold that will result in the least decrease in the minimum preparedness metric over all zones is dispatched, that is, the ambulance~$a$ with expected travel time less than a specified threshold that maximizes $\min\{\psi^{-a}_{\ell} : \ell \in \mathcal{L}\}$ is dispatched.
For ambulance reassignment decisions, an optimization problem is solved that minimizes the maximum travel time to increase the preparedness metric in all zones above a specified threshold.

\citet{liu:13} proposed the following modification of the preparedness metric $\psi_{\ell}$ of \citet{ande:07}, that takes into account that not all ambulances are equally likely to respond to an emergency in a zone.
Let $p^{a}_{\ell}$ denote the probability that ambulance~$a$ responds to an emergency in zone~$\ell$.
\citet{liu:13} used the M/M/N/loss-based model of \citet{budg:09} to compute $p^{a}_{\ell}$.
Then the modified preparedness metric is
\[
\psi^{p}_{\ell} \ \ \defi \ \ \frac{1}{\lambda_{\ell}} \sum_{a=1}^{A_{\ell}} \frac{\gamma^{a} p^{a}_{\ell}}{t^{a}_{\ell}}
\]
\citet{liu:13} used the modified preparedness metric in a chance constrained formulation for ambulance relocation.

Ambulance selection based on a preparedness metric was also considered by \citet{lees:11}.
It was demonstrated that the policy of \citet{ande:07} can result in worse performance than the closest available ambulance rule.
\citet{lees:11} used the same preparedness metric $\psi_{\ell}$ as \citet{ande:07}, but proposed two modifications to the use of the preparedness metric for ambulance selection decisions (only one priority level was considered).
The first modification selects the available ambulance that maximizes the minimum preparedness metric over all zones divided by the travel time from the ambulance to the emergency location, that is, the ambulance~$a$ that maximizes $\min\{\psi^{-a}_{\ell} : \ell \in \mathcal{L}\} / t^{a}_{e}$ ($\min\{\psi^{-a}_{\ell} : \ell \in \mathcal{L}\} / (1 + t^{a}_{e})$ in \citealt{lees:17}) is dispatched, where $t^{a}_{e}$ denotes the travel time of ambulance~$a$ from its current location to the location of the emergency~$e$.
The second modification considered in addition to the minimum preparedness metric over all zones $\min\{\psi^{-a}_{\ell} : \ell \in \mathcal{L}\}$ in the objective of the optimization problem, the following alternative aggregates of the preparedness metrics of different zones:
\begin{enumerate}
\item
The average preparedness metric over all zones $\overline{\psi}^{-a} \defi \sum_{\ell \in \mathcal{L}} \psi^{-a}_{\ell} / |\mathcal{L}|$.
\item
The average preparedness metric over all zones penalized by the inequality among zones, $\overline{\psi}^{-a} (1 - G^{-a})$, where $G^{-a}$ denotes the Gini index of the zone preparedness values $\{\psi^{-a}_{\ell} : \ell \in \mathcal{L}\}$ after removing ambulance~$a$.
\item
The average preparedness metric over all zones more heavily penalized by the inequality among zones, $\overline{\psi}^{-a} (1 - G^{-a}) / (1 + G^{-a})$.
\end{enumerate}
\citet{lees:17} pointed out the following shortcoming of the alternative aggregates above:
If two zones~$\ell, \ell' \in \mathcal{L}$ have the same travel times, $t^{a}_{\ell} = t^{a}_{\ell'}$ for all~$a$, and different demand rates, say $\lambda_{\ell} < \lambda_{\ell'}$, then $\psi_{\ell} > \psi_{\ell'}$, and thus the zone with less demand will have greater impact on the average preparedness metric $\overline{\psi} \defi \sum_{\ell \in \mathcal{L}} \psi_{\ell} / |\mathcal{L}|$.
To address the concern, \citet{lees:17} proposed an ambulance selection method that dispatches the ambulance~$a$ that minimizes the weighted average response time
\[
(1 + t^{a}_{e}) \left(\sum_{\ell \in \mathcal{L}} \lambda_{\ell} \left(1 + \min\left\{t^{a'}_{\ell} \; : \; a' \in \mathcal{A}_{\ell} \setminus \{a\}\right\}\right)\right)^{w}
\]
where $\mathcal{A}_{\ell} \subset \mathcal{A}$ denotes the set of (currently available) ambulances that can respond to an emergency in zone~$\ell$ and $w \ge 0$ is a user-chosen parameter that weighs preparedness relative to immediate travel time $t^{a}_{e}$.

\new{
\citet{carv:20,carv:23,carv:25} considered decision making in time periods.
More than one emergency may arrive in a time period, or more than one ambulance may finish service and become available in a time period.
Therefore, in each time period an optimization problem is considered that may select multiple ambulances for different emergencies and that may reassign (or reposition) multiple available ambulances.
The objective function includes both an immediate performance measure (the excess response time above a specified threshold response time) as well as a preparedness metric similar to the one proposed by \citet{lees:17}.
\citet{carv:23,carv:25} also considered a setting with multiple ambulance types, and proposed the following modification of the preparedness metric of \citet{lees:17}:
Let $A$ denote the set of ambulance types, and for each ambulance type $a \in A$, let $\lambda^{a}_{\ell}$ denote the demand for ambulances of type~$a$ in zone~$\ell$, and let $\mathcal{A}^{a}_{\ell}$ denote the set of (currently available) ambulances of type~$a$ that can respond to an emergency in zone~$\ell$.
Then the modified preparedness metric in the objective is
\[
\sum_{\ell \in \mathcal{L}} \sum_{a \in A} \lambda^{a}_{\ell} \; \min\left\{t^{a'}_{\ell} \; : \; a' \in \mathcal{A}^{a}_{\ell}\right\}
\]}

\citet{lees:12,lees:13} proposed a policy for assigning a newly available ambulance to emergency requests waiting in queue that takes into account both the travel times between the ambulance and requests in queue, as well as a centrality measure of each request in queue.
The idea is to give preference to a request in queue that is close to other requests in queue, so that if the ambulance does not have to take the patient of the first request to a hospital, then the ambulance will be close to other requests in queue after it finishes serving the first request.
Specifically, let $\mathcal{Q}$ denote the emergencies in queue, for each emergency~$e \in \mathcal{Q}$ let $t^{a}_{e}$ denote the travel time of the newly available ambulance~$a$ from its current location to the location of emergency~$e$, and for each pair of emergencies~$e,e' \in \mathcal{Q}$ let $t_{e,e'}$ denote the travel time between the locations of the emergencies.
\citet{lees:12} considered the following measures of the centrality~$c_{e}$ of an emergency~$e \in \mathcal{Q}$ relative to other emergencies~$e' \in \mathcal{Q}$ (\citealt{lees:13} restricted attention to the first of these centrality measures):
\begin{enumerate}
\item
The weighted degree $c_{e} = \sum_{e' \in \mathcal{Q} \setminus \{e\}} 1 / \left(1 + t_{e,e'}\right)$.
\item
The distance centrality $c_{e} = 1 / \left(1 + \sum_{e' \in \mathcal{Q} \setminus \{e\}} t_{e,e'}\right)$.
\item
The betweenness centrality $c_{e} = \sum_{e',e'' \in \mathcal{Q} \setminus \{e\}} \sigma_{e',e,e''} / \sigma_{e',e''}$, where $\sigma_{e',e''}$ is the number of shortest paths (with equal travel times) between $e',e'' \in \mathcal{Q}$, and $\sigma_{e',e,e''}$ is the number of these shortest paths that go through the location of~$e$, that is, $\sigma_{e',e,e''} / \sigma_{e',e''}$ is the fraction of shortest paths between $e',e'' \in \mathcal{Q}$ that go through~$e$.
\end{enumerate}
For any chosen centrality measure $c_{e}$ and (exponent) weight~$w > 0$ (for example, $w =$ the probability that the ambulance does not transport the patient to a hospital), if an ambulance~$a$ becomes available and there are emergencies in queue, then the ambulance is dispatched to emergency $e^* \in \arg\max\left\{c_{e}^{w} / (1 + t^{a}_{e}) : e \in \mathcal{Q}\right\}$.
In addition, \citet{lees:14} took into account both available and busy ambulances when making ambulance dispatch decisions.
(In such a setting, there is no distinction between the ambulance selection problem and the ambulance reassignment problem.)
Specifically, let $\mathcal{A}$ denote the set of all ambulances, both available and busy.
If ambulance~$a$ is available, then as before $t^{a}_{e}$ denotes the travel time of ambulance~$a$ from its current location to the location of emergency~$e$, and if ambulance~$a$ is busy, then $t^{a}_{e}$ denotes the forecasted remaining time for ambulance~$a$ to complete its current service plus the travel time of ambulance~$a$ from where it will complete its current service to the location of emergency~$e$.
A centrality measure $c_{e}$ is calculated for each emergency~$e \in \mathcal{Q}$ as described above.
Then the assignment problem
\[
\max_{x \in \{0,1\}^{\mathcal{A} \times \mathcal{Q}}} \left\{\sum_{a \in \mathcal{A}} \sum_{e \in \mathcal{Q}} \frac{c_{e}^{w}}{(1 + t^{a}_{e})} x_{a,e} \; : \; \sum_{a \in \mathcal{A}} x_{a,e} \le 1 \; \forall \; e \in \mathcal{Q}, \; \sum_{e \in \mathcal{Q}} x_{a,e} \le 1 \; \forall \; a \in \mathcal{A}\right\}
\]
is solved, producing optimal solution $x^*$.
For each available ambulance~$a$ and emergency~$e$, if $x^*_{a,e} = 1$ then ambulance~$a$ is dispatched to emergency~$e$, while for busy ambulances, decisions are postponed.

Ambulance selection policies based on approximations of consequences were also proposed in \citet{jagt:17a,jagt:17b}.
One selection policy uses a simplified discrete-time Markov decision process, based on the assumption that ambulance combined travel and service times are geometrically distributed with mean independent of the ambulance location or emergency location.
To keep the state space relatively small, the state keeps track only of the busy/idle status of each ambulance, and the location of a newly arrived emergency.
The other selection policy is based on a preparedness metric similar to the objective value of the maximum expected covering location problem (MEXCLP) of \citet{dask:82,dask:83}, as follows.
Let $\mathcal{A}$ denote the set of currently available ambulances, and for each zone~$\ell \in \mathcal{L}$, let $\mathcal{A}_{\ell} \subset \mathcal{A}$ denote the set of (currently available) ambulances that can respond to an emergency in zone~$\ell$ in less time than a specified threshold~$T$.
Let $\lambda_{\ell}$ denote the demand rate for zone~$\ell$, and let $q \in [0,1]$ denote the ``busy fraction'' of all ambulances, i.e., the fraction of time that ambulances are busy.
If an emergency occurs in zone~$\ell'$ and $\mathcal{A}_{\ell'} \neq \varnothing$, then an ambulance in $\mathcal{A}_{\ell'}$ that minimizes the reduction in the preparedness metric is chosen, i.e.,
\[
a^* \ \ \in \ \ \arg\min\left\{\sum_{\{\ell \in \mathcal{L} \, : \, a \in \mathcal{A}_{\ell}\}} \lambda_{\ell} (1 - q) q^{|\mathcal{A}_{\ell}| - 1} \; : \; a \in \mathcal{A}_{\ell'}\right\}
\]
is chosen.
If an emergency occurs in zone~$\ell'$ and $\mathcal{A}_{\ell'} = \varnothing$, then
\[
a^* \ \ \in \ \ \arg\min\left\{\sum_{\{\ell \in \mathcal{L} \, : \, a \in \mathcal{A}_{\ell}\}} \lambda_{\ell} (1 - q) q^{|\mathcal{A}_{\ell}| - 1} \; : \; a \in \mathcal{A}\right\}
\]
is chosen.
Simulation results showed that the policy has a much lower fraction of late arrivals than the closest available ambulance rule, but that the policy also has a much greater mean response time than the closest available ambulance rule.

\new{
The hypercube queueing model of \citet{lars:74} used a continuous-time Markov process, similar to the model in this paper.
\citet{lars:74} considered a setting with $N$ distinct servers, multiple demand locations, and a given fixed preference policy, that is, a policy that specifies for every demand point a preference list of all the servers from most preferred to least preferred (with ties allowed) independent of the state of the process.
When a call arrives from a demand point, the most preferred available server in the preference list for that demand point (one of the most preferred available servers in case of ties) is dispatched to serve the call.
An algorithm was proposed to compute the transition rates.
The algorithm exploits the similarity of the most preferred available servers for adjacent states of the Markov process, that is, states that differ in the availability of only one server, to reduce the effort to compute the transition rates.
After the transition rates have been computed, a system of $2^N$ linear equations can be solved to compute the stationary probabilities, and then various long-run average performance metrics can be computed.
Unlike the model of \citet{lars:74}, our model makes provision for different emergency types, ambulance types, and ambulance stations, but Markov chains are separated by ambulance station.}

\subsubsection{Different performance metrics}

Many EMSs as well as academic papers use response time as a performance metric for ambulance operations.
Usually, summary statistics of response time data, such as average response time and specific quantiles of the response time empirical distribution, are considered.
For example, it is often checked whether the $0.8$ or $0.9$ empirical response time quantiles are less than specified threshold values \citep{hend:04,rest:09,maso:13}.
\citet{erku:08,band:12,knig:12b,mayo:13,band:14} considered models that maximize patient survival probability, as opposed to models that maximize the probability that the response time is less than a threshold.

The relationship between EMS response time and the probability of survival to hospital admission or the probability of survival to hospital discharge has been studied, for example in \citet{cret:79,blac:91,lars:93,vale:97,stie:99,vale:00,waal:01,pell:01,pons:02,dema:03,pons:05,stie:08,blac:09,okee:11,blan:12,weis:13,wild:13}.
For out of hospital cardiac arrest, the results support the hypothesis that the probability of survival decreases as the elapsed time from collapse to cardiopulmonary resuscitation (CPR), or to first defibrillation, or to initiation of advanced cardiac life support, increases \citep{lars:93,vale:97,stie:99,vale:00,pell:01,waal:01,dema:03}.
For other types of emergencies, studies often found no statistically significant relationship between response time and survival probability \citep{blac:09,blan:12,weis:13}.
Thus, the relationship between response time and survival probability depends on the emergency type, and except possibly for cardiac arrest, much more research is needed to accurately quantify the relationship.

It has also been pointed out that in addition to response time, the capabilities of the ambulance and personnel can affect the survival probability of the patient, depending on the emergency type.
The dispatcher may choose among different ambulance types for dispatch, such as a BLS or an ALS ambulance, which may contain crew members with different qualifications, such as EMTs, AEMTs, paramedics, or physicians, and the best choice may depend on the type of emergency. 
Our models allow the performance metric to depend on the emergency type, equipment/crew combination, and response time, and thus make provision for various predicted patient outcome measures.
For most emergency types, accurate prediction models are not available, and therefore our numerical results serve as examples only, and not as guidelines for dispatching decisions.

\if{

\section{Notation and terminology}

A crucial step, that applies for all the policies proposed in this paper, is to define the state vector of the system at each time instant.
This is the minimal information necessary to take decisions.
For all policies, the state vector will store, at any time $t$:
\begin{enumerate}
\item[1)]
for each ambulance $j=1,\ldots,$ NbAmbulances, a location $\ell_{f}(j)$ and a time $t_{f}(j)$ with the following meanings.
Two situations can happen:
1.1) an ambulance $j$ is in service at $t$ or 1.2) it is available at $t$.
\par 1.1) If the ambulance is in service at $t$ then $t_{f}(j)$ is the next instant the ambulance will be available (free) again for dispatch and $\ell_{f}(j)$ will be its location at that instant.
Therefore, $t_{f}(j)>t$ (the strict inequality holds because the ambulance is not currently available) is a future time instant and $\ell_{f}(j)$ is a future ambulance location.
\par 1.2) If the ambulance is available at $t$ then $t_{f}(j)$ is the last (past) time instant the ambulance became available, i.e., the time it completed its last service and  $\ell_{f}(j)$ was the corresponding (past) location of the ambulance when this service was completed (either the location of a call, of a hospital, or of a cleaning base).
\item[2)]
For each ambulance $j=1,\ldots,$ NbAmbulances, a location $\ell_{b}(j)$ and a time $t_{b}(j)$ with the following meanings.
We again have two possibilities.
Either the ambulance is at a base at $t$ and, in this case, $\ell_{b}(j)$ is the location of that base and $t_{b}(j)$ is the last time instant it arrived at this base (therefore $t_{b}(j) \leq t$).
If the ambulance is not at a base, then $\ell_{b}(j)$ is the location of the next base it will go and $t_{b}(j)$ is the instant the ambulance will arrive at that base (we therefore have $t_{b}(j) > t$).
\end{enumerate}

The $i$th call over the planning horizon is given by:
\begin{itemize}
\item
the instant $t_{c}(i)$ of the call;
\item
the location $\ell_{c}(i)$ of the call (given by a pair latitude/longitude);
\item
$\mbox{Type}_{c}(i)$: the type of call $i$ (a positive integer)\footnote{Usually, the types of calls are grouped into groups of calls of different priorities and the priority of the call is enough to decide which ambulance types are possible for that call.};
\item
$\mbox{TimeOnScene}_{c}(i)$: time spent by the ambulance on the scene of the call;
\item
$h_{c}(i)$: the index of the hospital the patient of call $i$ is sent to;\footnote{When a hospital is assigned to a call already at the instant of the call, it is usually the closest hospital to the scene of the call that can attend the corresponding emergency type.
Other criteria for choosing the hospitals, taking into account both the time to go to the hospital and bed availabilities in hospitals will be taken into account in Section~\ref{sec:optqueue}.}
\item
$\mbox{HospitalNeeded}_{c}(i)$: is 1 if the corresponding patient is transported to a hospital and $0$ otherwise;
\item
$\mbox{TimeAtHospital}_{c}(i)$: time the ambulance waits the bed to be freed at hospital after transporting the patient.
\item
$\mbox{CleaningNeeded}_{c}(i)$: is 1 if the ambulance attending that call needs cleaning after serving that call and $0$ otherwise;
\item
$\mbox{CleaningTime}_{c}(i)$: time required to clean the ambulance after serving call $i$.
\end{itemize}

For each ambulance $j$, we define
\begin{itemize}
\item
its type $\mbox{Type}(j)$ which is a positive integer.
\end{itemize}

Each parking base $i$ has a location Base($i$) while cleaningBase($i$).location is the location of cleaning base $i$.
The location (given by a pair (latitude, longitude)) of the closest parking base to cleaning base $i$ is cleaningBase(i).base.

We will also denote by
\begin{enumerate}
\item
hospitals($i$).location the location (given by a pair (latitude/longitude)) of hospital $i$,
\item
hospitals($i$).base the base which is the closest to hospital $i$.
\item
hospitals($i$).cleaningBase the index of the cleaning base which is the closest to hospital $i$.
\end{enumerate}

}\fi

\subsection{Contributions}
\label{sec:contributions}

The main contribution of this paper is to propose $2$~new policies for the operation of an ambulance fleet under uncertainty.
The policies are based on a new preparedness metric described in Section~\ref{sec:preparedness}.
Whereas the preparedness metrics reviewed in the previous section take into account immediate considerations such as the currently available ambulances, the new preparedness metric models how the number of available ambulances at a station vary over time.
Also, whereas previous preparedness metrics are either linear in demand $\lambda_{\ell}$ \citep{lees:17,jagt:17a,jagt:17b}, or linear in a measure of supply relative to demand \citep{ande:07,lees:11}, the new preparedness metric allows performance to be nonlinear in supply relative to demand, as is typical in queueing systems.
In addition, the new preparedness metric makes provision for different types of emergencies and different types of ambulances and crews.
The policy for ambulance selection is described in Section~\ref{sec:selection}, and the policy for ambulance reassignment is described in Section~\ref{sec:reassignment}.
The policies are nonanticipative, are computed quickly in real time, and numerical results show that they provide good performance for emergency medical services.
Specifically, numerical results based on real data for a large EMS are presented in Section~\ref{sec:numsec}, where we show that in most cases the proposed method outperforms $9$~policies from the literature.

\section{A New Preparedness Metric}
\label{sec:preparedness}

In this section, we propose a continuous-time Markov chain model to measure the preparedness of an ambulance station with a given set of ambulances at the station.
The model is a local approximation in the sense that it considers only a single station, and it considers only the currently available ambulances at the station (thus it is local in both space and time).
The benefit of this local approximation is that it is tractable while providing an accurate metric of the ability to dispatch ambulances to emergencies in the zone of the station without delay as long as ambulance supply is sufficient relative to demand.
This approximation may be inaccurate if demand is large relative to supply, so that stations often run out of ambulances and ambulances often have to travel long distances to take care of emergencies.
Good dispatch policies attempt to prevent ambulances having to travel long distances to emergencies, and thus as long as the EMS has sufficient supply and good dispatch policies are used, the local approximation describes how ambulance operations take place most of the time.

The model considers a single station at a time, and therefore the notation does not indicate the station being considered.
For each emergency type $c \in C$, let $\lambda(c)$ denote the arrival rate of type~$c$ emergencies to the zone of the station, and let $A(c)$ denote the set of ambulance types that can serve a type~$c$ emergency.
Conversely, for each ambulance type $a \in A$, let $C(a) \defi \{c \in C \, : \, a \in A(c)\}$ denote the set of emergency types that can be served by ambulance type~$a$.
We assume that associated with each emergency type $c \in C$ there is a total order $\succ_{c}$ on $A(c)$ that governs the type of ambulance dispatched in the Markov chain model (not necessarily in the policy proposed in Section~\ref{sec:selection}) to serve a type~$c$ emergency, as follows.
If a type~$c$ emergency arrives, and an ambulance of type $a \in A(c)$ is available, then such an ambulance is dispatched to the emergency in the Markov chain model if and only if no ambulance of type $a' \in A(c)$ such that $a' \succ_{c} a$ is available, that is, the most preferred ambulance type for the emergency type that is available is dispatched in the Markov chain model.
Note that $a' \succ_{c} a$ does not mean that type~$a'$ ambulances are more capable than type~$a$ ambulances.
For example, for a minor emergency type~$c$ it may hold that BLS $\succ_{c}$ ALS.
If a type~$a$ ambulance is dispatched to a type~$c$ emergency in the zone of the station, then the ambulance becomes available again after an exponentially distributed amount of time with mean $1 / \mu(a,c)$ (in the Markov chain model).

For each ambulance type $a \in A$, let $m(a)$ denote the number of type~$a$ ambulances at the station, and let $m \defi (m(a), a \in A)$.
The state of the Markov chain is specified by the number of ambulances of each type that is busy serving emergencies of each type.
For each emergency type $c \in C$ and ambulance type $a \in A(c)$, let $x_{a,c}$ denote the number of type~$a$ ambulances currently busy serving type~$c$ emergencies.
It is assumed that each ambulance not busy serving an emergency is available, and thus the number of available type~$a$ ambulances is given by $m(a) - \sum_{c \in C(a)} x_{a,c}$.
Let $x \defi (x_{a,c}, a \in A, c \in C(a))$ denote the state of the Markov chain.

Let $q^{m}_{x,x'}$ denote the transition rate from state~$x$ to state~$x'$.
There are two types of state transitions: transitions caused by emergency arrivals and subsequent ambulance dispatches, and transitions caused by ambulances completing tasks and becoming available.
Consider any current state~$x$ and any emergency type~$c$ such that there is an ambulance available to serve an emergency of that type, that is, $\sum_{a \in A(c)} \left[m(a) - \sum_{c' \in C(a)} x_{a,c'}\right] > 0$.
Given state~$x$, let
\[
\hat{a} \ \ \in \ \ \left\{a \in A(c) \; : \; m(a) > \sum_{c' \in C(a)} x_{a,c'}, \; \nexists \, a' \in A(c) \mbox{ s.t.\ } m(a') > \sum_{c' \in C(a')} x_{a',c'}, \; a' \succ_{c} a\right\}
\]
denote the available ambulance type that is preferred in state~$x$ for type~$c$ emergencies.
Note that $\hat{a}$ depends on~$x$ and~$c$, but the notation does not show the dependence.
Then $q^{m}_{x,x'} = \lambda(c)$, where $x'$ is determined as follows:
\begin{eqnarray*}
x'_{\hat{a},c} & = & x_{\hat{a},c} + 1 \\
x'_{a',c'} & = & x_{a',c'} \qquad \forall \ (a',c') \neq (\hat{a},c).
\end{eqnarray*}
In addition, for each emergency type $c \in C$ and ambulance type $a \in A(c)$ such that $x_{a,c} > 0$, $q^{m}_{x,x'} = \mu(a,c) x_{a,c}$, where $x'$ is determined as follows:
\begin{eqnarray*}
x'_{a,c} & = & x_{a,c} - 1 \\
x'_{a',c'} & = & x_{a',c'} \qquad \forall \ (a',c') \neq (a,c).
\end{eqnarray*}
Let $q^{m}_{x,x} \defi - \sum_{x'} q^{m}_{x,x'}$.
Note that all states communicate, for example, any state~$x'$ can be reached from state~$x = 0$ in a finite number of transitions and vice versa.
Thus, the Markov chain has a unique stationary distribution~$\nu^{m} = (\nu^{m}_{x})$ given by
\begin{eqnarray}
\label{eqn:balance equations}
\sum_{x} \nu^{m}_{x} q^{m}_{x,x'}  & = & 0 \qquad \forall \ x' \\
\sum_{x} \nu^{m}_{x} & = & 1.
\label{eqn:probability scaling}
\end{eqnarray}

Given stationary distribution~$\nu^{m} = (\nu^{m}_{x})$, one can compute steady-state performance measures such as the following.
Let $\phi(c)$ denote the cost/penalty if there is no ambulance available to serve a type~$c$ emergency.
Then the cost rate $\psi_{x}$ while in state~$x$ is given by
\[
\psi^{m}_{x} \ \ = \ \ \sum_{c \in C} \mathds{1}\left\{\sum_{a \in A(c)} \left[m(a) - \sum_{c' \in C(a)} x_{a,c'}\right] = 0\right\} \lambda(c) \phi(c)
\]
and the steady state cost rate $\bar{\psi}^{m}$ given ambulance supply~$m$ is $\bar{\psi}^{m} = \sum_{x} \nu^{m}_{x} \psi^{m}_{x}$.
More generally, any state-dependent performance measure $\psi^{m}_{x}$ can be specified and converted to a steady-state performance measure $\bar{\psi}^{m}$.

For each station~$b$, $\bar{\psi}^{m}_{b}$ can be computed in advance for each $m \defi (m(a), a \in A)$.
Then $\bar{\psi}^{m}_{b}$ can be used as a ``preparedness'' metric for ambulance supply~$m$ at station~$b$.
For example, when an ambulance becomes available after completing a task, it can be sent to the station~$b$ where it will improve $\bar{\psi}^{m}_{b}$ the most.
More specifically, for each station~$b$, let $m_{b}$ denote the current ambulance supply at station~$b$,
and let $m^{+}_{b}$ denote the ambulance supply at station~$b$ if the newly available ambulance would be added to the ambulance supply at station~$b$.
Then send the newly available ambulance to the station $b^* \in \arg\max\left\{\bar{\psi}^{m_{b}}_{b} - \bar{\psi}^{m^{+}_{b}}_{b} \, : \, b \in B\right\}$.

\section{A New Ambulance Selection Policy}
\label{sec:selection}

In this section, we propose a new ambulance selection policy.
Suppose that a type~$c$ emergency arrives at location~$\ell$.
Let $i_{0}$ be the index of the newly arrived emergency.
For any station~$b$, let $\mathcal{A}_{b}$ denote the set of ambulances currently at station~$b$ or en route to station~$b$, and let $\mathcal{A}'$ denote the set of currently on-task ambulances.
Thus the set of ambulances is $\mathcal{A} = \cup_{b \in B} \mathcal{A}_{b} \cup \mathcal{A}'$.
For each ambulance $a \in \mathcal{A}$, let $t(a) \in A$ denote the type of the ambulance.
For each station~$b$, let $m_{b} = (m_{b}(\tilde{a}), \tilde{a} \in A)$ denote the current ambulance supply at station~$b$, with the understanding that any ambulance currently en route to a station is included in the ambulance supply of that station.
Let $\mathcal{Q}$ denote the set of emergencies currently in queue.
For each emergency $i \in \mathcal{Q} \cup \{i_{0}\}$ and each ambulance $a \in \mathcal{A}$, let $r(a,i)$ denote the immediate impact on the performance measure if ambulance~$a$ is dispatched to emergency~$i$; for example, $r(a,i)$ may represent the response time for emergency~$i$ if ambulance~$a$ is dispatched to the emergency.
Note that if ambulance~$a$ is currently on-task, then the calculation of $r(a,i)$ is based on ambulance~$a$ first completing its current task and then traveling to emergency~$i$.

For each emergency $i \in \mathcal{Q} \cup \{i_{0}\}$ and each ambulance $a \in \mathcal{A}$, let $x(a,i)$ denote a decision variable that is~$1$ if ambulance~$a$ is dispatched to emergency~$i$, and is~$0$ otherwise.
For each on-task ambulance $a \in \mathcal{A}'$ and station~$b \in B$, let $y(a,b)$ denote a decision variable that is~$1$ if ambulance~$a$ is sent to station~$b$, and is~$0$ otherwise.
Let $x \defi (x(a,i), a \in \mathcal{A}, i \in \mathcal{Q} \cup \{i_{0}\})$ and $y \defi (y(a,b), a \in \mathcal{A}', b \in B)$.
For any ambulance type~$\tilde{a} \in A$ and station~$b \in B$, let $m^{+}_{b}(\tilde{a},x,y) \defi m_{b}(\tilde{a}) - \sum_{a \in \mathcal{A}_{b}} \mathds{1}\left\{t(a) = \tilde{a}\right\} \sum_{i \in \mathcal{Q} \cup \{i_{0}\}} x(a,i) + \sum_{a \in \mathcal{A}'} \mathds{1}\left\{t(a) = \tilde{a}\right\} y(a,b)$, and let $m^{+}_{b}(x,y) \defi (m^{+}_{b}(\tilde{a},x,y), \tilde{a} \in A)$.

Let $\Gamma$ be a parameter that weighs the effects of preparedness (uncertain future response times) relative to current response times.
For each emergency $i \in \mathcal{Q} \cup \{i_{0}\}$, let $\gamma(i)$ be a parameter that weighs the effect of queueing emergency~$i$; $\gamma(i)$ should depend on the type of emergency~$i$.
Version~1 of the ambulance selection policy solves the following optimization problem:
\begin{eqnarray}
\min_{x,y} & & \sum_{a \in \mathcal{A}} \sum_{i \in \mathcal{Q} \cup \{i_{0}\}} r(a,i) x(a,i) + \sum_{i \in \mathcal{Q} \cup \{i_{0}\}} \gamma(i) \left[1 - \sum_{a \in \mathcal{A}} x(a,i)\right] + \Gamma \sum_{b \in \mathcal{B}} \bar{\psi}^{m^{+}_{b}(x,y)}_{b}
\label{eqn:nonlinear selection policy objective} \\
\mbox{s.t.} & & \sum_{i \in \mathcal{Q} \cup \{i_{0}\}} x(a,i) \ \ \le \ \ 1 \qquad \forall \ a \in \cup_{b \in B} \mathcal{A}_{b}
\label{eqn:selection each ambulance at most one dispatch} \\
& & \sum_{i \in \mathcal{Q} \cup \{i_{0}\}} x(a,i) + \sum_{b \in B} y(a,b) \ \ = \ \ 1 \qquad \forall \ a \in \mathcal{A}'
\label{eqn:selection each on-task ambulance exactly one dispatch} \\
& & \sum_{a \in \mathcal{A}} x(a,i) \ \ \le \ \ 1 \qquad \forall \ i \in \mathcal{Q} \cup \{i_{0}\}.
\label{eqn:selection each emergency at most one dispatch}
\end{eqnarray}

The first term in the objective function~\eqref{eqn:nonlinear selection policy objective} represents the immediate impact on the performance measure of dispatching ambulances to emergencies in queue, including the emergency call~$i_{0}$ that has just arrived.
The second term in the objective function represents the impact on the performance measure of not yet dispatching an ambulance to an emergency in queue, that is, an emergency in queue waiting longer before dispatching an ambulance to the emergency.
The third term in the objective function represents the impact of the current dispatch decisions on the preparedness metric, that is, the impact of the current dispatch decisions on the expected future performance measure.
Constraint~\eqref{eqn:selection each ambulance at most one dispatch} restricts each ambulance at a station to be immediately dispatched to at most one emergency.
Constraint~\eqref{eqn:selection each on-task ambulance exactly one dispatch} requires that each ambulance not at a station be immediately dispatched either to an emergency or to a station.
Constraint~\eqref{eqn:selection each emergency at most one dispatch} requires at most one ambulance to be dispatched to an emergency.
The objective function~\eqref{eqn:nonlinear selection policy objective} is not linear in $(x,y)$.
Version~2 of the policy solves a linear optimization problem.
For each ambulance type~$\tilde{a} \in A$, let $e(\tilde{a})$ denote the unit vector with the component for $\tilde{a}$ equal to~$1$ and the components for all $\tilde{a}' \in A \setminus \{\tilde{a}\}$ equal to~$0$.
For each ambulance type~$\tilde{a} \in A$ and station~$b \in B$, let $s^{+}(\tilde{a},b) \defi \bar{\psi}^{m_{b} + e(\tilde{a})}_{b} - \bar{\psi}^{m_{b}}_{b}$ and $s^{-}(\tilde{a},b) \defi \bar{\psi}^{m_{b} - e(\tilde{a})}_{b} - \bar{\psi}^{m_{b}}_{b}$.
Then Version~2 of the ambulance selection policy solves the following linearized optimization problem:
\begin{eqnarray}
\min_{x,y} & & \sum_{a \in \mathcal{A}} \sum_{i \in \mathcal{Q} \cup \{i_{0}\}} r(a,i) x(a,i) + \sum_{i \in \mathcal{Q} \cup \{i_{0}\}} \gamma(i) \left[1 - \sum_{a \in \mathcal{A}} x(a,i)\right] \nonumber \\
& & \ {} + \Gamma \sum_{b \in \mathcal{B}} \sum_{a \in \mathcal{A}_{b}} \sum_{i \in \mathcal{Q} \cup \{i_{0}\}} s^{-}(t(a),b) x(a,i) + \Gamma \sum_{a \in \mathcal{A}'} \sum_{b \in \mathcal{B}} s^{+}(t(a),b) y(a,b)
\label{eqn:linear selection policy objective} \\
\mbox{s.t.} & & \sum_{i \in \mathcal{Q} \cup \{i_{0}\}} x(a,i) \ \ \le \ \ 1 \qquad \forall \ a \in \cup_{b \in B} \mathcal{A}_{b}
\label{eqn:linear selection each ambulance at most one dispatch} \\
& & \sum_{i \in \mathcal{Q} \cup \{i_{0}\}} x(a,i) + \sum_{b \in B} y(a,b) \ \ = \ \ 1 \qquad \forall \ a \in \mathcal{A}'
\label{eqn:linear selection each on-task ambulance exactly one dispatch} \\
& & \sum_{a \in \mathcal{A}} x(a,i) \ \ \le \ \ 1 \qquad \forall \ i \in \mathcal{Q} \cup \{i_{0}\}.
\label{eqn:linear selection each emergency at most one dispatch}
\end{eqnarray}

The first two terms in the objective function~\eqref{eqn:linear selection policy objective} are the same as in~\eqref{eqn:nonlinear selection policy objective}.
The third and fourth terms in the objective function represent the linearized impact of the current dispatch decisions on the expected future performance measure, calculated through the differences made by each dispatch decision on the preparedness metric.
Constraints~\eqref{eqn:linear selection each ambulance at most one dispatch}--\eqref{eqn:linear selection each emergency at most one dispatch} are the same as constraints~\eqref{eqn:selection each ambulance at most one dispatch}--\eqref{eqn:selection each emergency at most one dispatch}.

Let $(x^*,y^*)$ denote an optimal solution of Version~1 or Version~2 of the ambulance selection policy.
With both versions of the policy, the ambulance selection decision dispatches ambulance~$a$ to newly arrived emergency~$i_{0}$ if and only if $x^*(a,i_{0}) = 1$.

\section{A New Ambulance Reassignment Policy}
\label{sec:reassignment}

In this section, we propose a new ambulance reassignment policy.
Suppose that ambulance~$a_{0} \in \mathcal{A}$ becomes available after completing a task.
The rest of the notation is the same as in Section~\ref{sec:selection}.
Version~1 of the ambulance reassignment policy solves the following optimization problem:
\begin{eqnarray}
\min_{x,y} & & \sum_{a \in \mathcal{A}} \sum_{i \in \mathcal{Q}} r(a,i) x(a,i) + \sum_{i \in \mathcal{Q}} \gamma(i) \left[1 - \sum_{a \in \mathcal{A}} x(a,i)\right] + \Gamma \sum_{b \in \mathcal{B}} \bar{\psi}^{m^{+}_{b}(x,y)}_{b}
\label{eqn:nonlinear reassignment policy objective} \\
\mbox{s.t.} & & \sum_{i \in \mathcal{Q}} x(a,i) \ \ \le \ \ 1 \qquad \forall \ a \in \cup_{b \in B} \mathcal{A}_{b}
\label{eqn:reassignment each ambulance at most one dispatch} \\
& & \sum_{i \in \mathcal{Q}} x(a,i) + \sum_{b \in B} y(a,b) \ \ = \ \ 1 \qquad \forall \ a \in \mathcal{A}' \cup \{a_{0}\}
\label{eqn:reassignment each on-task ambulance exactly one dispatch} \\
& & \sum_{a \in \mathcal{A}} x(a,i) \ \ \le \ \ 1 \qquad \forall \ i \in \mathcal{Q}
\label{eqn:reassignment each emergency at most one dispatch}
\end{eqnarray}

Problem~\eqref{eqn:nonlinear reassignment policy objective}--\eqref{eqn:reassignment each emergency at most one dispatch} is the same as problem~\eqref{eqn:nonlinear selection policy objective}--\eqref{eqn:selection each emergency at most one dispatch} except for the following:
In ambulance selection problem~\eqref{eqn:nonlinear selection policy objective}--\eqref{eqn:selection each emergency at most one dispatch} there is a newly arrived emergency~$i_{0}$, and in ambulance reassignment problem~\eqref{eqn:nonlinear reassignment policy objective}--\eqref{eqn:reassignment each emergency at most one dispatch} there is not, whereas in ambulance reassignment problem~\eqref{eqn:nonlinear reassignment policy objective}--\eqref{eqn:reassignment each emergency at most one dispatch} there is a newly available ambulance~$a_{0}$, and in ambulance selection problem~\eqref{eqn:nonlinear selection policy objective}--\eqref{eqn:selection each emergency at most one dispatch} there is not.

As for~\eqref{eqn:nonlinear selection policy objective}, the objective function~\eqref{eqn:nonlinear reassignment policy objective} is not linear in $(x,y)$.
Therefore Version~2 of the ambulance reassignment policy solves the following linearized optimization problem:
\begin{eqnarray}
\min_{x,y} & & \sum_{a \in \mathcal{A}} \sum_{i \in \mathcal{Q}} r(a,i) x(a,i) + \sum_{i \in \mathcal{Q}} \gamma(i) \left[1 - \sum_{a \in \mathcal{A}} x(a,i)\right] \nonumber \\
& & \ {} + \Gamma \sum_{b \in \mathcal{B}} \sum_{a \in \mathcal{A}_{b}} \sum_{i \in \mathcal{Q}} s^{-}(t(a),b) x(a,i) + \Gamma \sum_{a \in \mathcal{A}' \cup \{a_{0}\}} \sum_{b \in \mathcal{B}} s^{+}(t(a),b) y(a,b)
\label{eqn:linear reassignment policy objective} \\
\mbox{s.t.} & & \sum_{i \in \mathcal{Q}} x(a,i) \ \ \le \ \ 1 \qquad \forall \ a \in \cup_{b \in B} \mathcal{A}_{b} \\
& & \sum_{i \in \mathcal{Q}} x(a,i) + \sum_{b \in B} y(a,b) \ \ = \ \ 1 \qquad \forall \ a \in \mathcal{A}' \cup \{a_{0}\} \\
& & \sum_{a \in \mathcal{A}} x(a,i) \ \ \le \ \ 1 \qquad \forall \ i \in \mathcal{Q}
\label{eqn:linear reassignment policy at most one ambulance}
\end{eqnarray}
Let $(x^*,y^*)$ denote an optimal solution of Version~1 or Version~2 of the ambulance reassignment policy.
With both versions of the policy, the ambulance reassignment decision dispatches ambulance~$a_{0}$ to emergency~$i \in \mathcal{Q}$ if and only if $x^*(a_{0},i) = 1$, and dispatches ambulance~$a_{0}$ to station~$b \in B$ if and only if $y^*(a_{0},b) = 1$.

\section{Numerical Results}
\label{sec:numsec}

The system of equations~\eqref{eqn:balance equations}--\eqref{eqn:probability scaling} may be quite large, and thus care should be taken to solve the system.
Note that the system of equations~\eqref{eqn:balance equations}--\eqref{eqn:probability scaling} has one more equation than unknown, and the equation $\sum_{x} \nu^{m}_{x} q^{m}_{x,x'} = 0$ for one $x'$, such as $x' = 0$, can be dropped.
Let $Q^{m} \defi (q^{m}_{x,x'})$ denote the transition rate matrix, let $\mathbf{0}$ denote the vector $(0,\ldots,0)$, and let $\mathbf{1}$ denote the vector $(1,\ldots,1)$.
Then the system~\eqref{eqn:balance equations}--\eqref{eqn:probability scaling} can be written as ${Q^{m}}^{\top} \nu^{m} = \mathbf{0}, \mathbf{1}^{\top} \nu^{m} = 1$.
Let $D^{m} \defi (\mathbf{1}, (q^{m}_{x,x'}, x' \neq 0))$ denote the reduced matrix.
Then the system~\eqref{eqn:balance equations}--\eqref{eqn:probability scaling} is equivalent to ${D^{m}}^{\top} \nu^{m} = (1,0,\ldots,0)$.
Note that the transition rate matrix~$Q^{m}$ is very sparse, and thus~$D^{m}$ is sparse.
We used the following two methods to solve the system~\eqref{eqn:balance equations}--\eqref{eqn:probability scaling}:
\begin{enumerate}
\item
The Generalized Minimal Residual method (GMRES) of \citet{saad:86} was used to solve the system ${D^{m}}^{\top} \nu^{m} = (1,0,\ldots,0)$ with asymmetric~$D^{m}$.
\item
The conjugate gradient method was used to solve the system 
$$D^{m} {D^{m}}^{\top} \nu^{m} = D^{m} (1,0,\ldots,0) = \mathbf{1}.
$$
The matrix $D^{m} {D^{m}}^{\top}$ has large diagonal entries, and thus the diagonal of $D^{m} {D^{m}}^{\top}$ was used as a preconditioner.
Also, the matrix $D^{m} {D^{m}}^{\top}$ was still quite sparse (but not as sparse as $D^{m}$).
\end{enumerate}
Both these methods exploit the sparsity of~$D^{m}$.

We compared the performance of different policies using data of Rio de Janeiro EMS.
The data include the history of emergency calls for the period January 2016--February 2018, the locations of ambulance stations and hospitals, and the set of ambulances.
We consider $2$~types of ambulances: Basic Life Support (BLS) and Advanced Life Support (ALS), and $4$~types of emergencies, numbered 1, 2, 3 and 4.
Any ambulance may serve any call, but every emergency has a time urgency (high or low) and an ambulance preference, as follows: 
\begin{itemize}
\item
type 1 call: a high-priority emergency that should preferably be served by an ALS ambulance;
\item
type 2 call: a low-priority emergency that should preferably be served by an ALS ambulance;
\item
type 3 call: a high-priority emergency that any ambulance can serve;
\item
type 4 call: a low-priority emergency that any ambulance can serve.
\end{itemize}
To model these preferences between emergency types and ambulance types, we use function \newline
cost\_allocation\_ambulance given by
\begin{equation}
\label{costalloc}
\mbox{cost\_allocation\_ambulance}(a,c,t) \ \ = \ \ \mbox{{\tt{penalization}}}(t,c) + M_{ac}
\end{equation}
which specifies the cost of allocating an ambulance of type~$a$ to an emergency of type~$c$ with response time~$t$ (the time elapsed from the instant the emergency call is received and the instant the ambulance arrives on the scene of the emergency).
\begin{itemize}
\item
In \eqref{costalloc}, {\tt{penalization}}$(t,c)$ is the penalized response time given by
\begin{equation}
\label{penbtct}
\mbox{{\tt{penalization}}}(t,c) \ \ = \ \ \theta_{c} t
\end{equation}
where $\theta_{c}$ is a coefficient depending on emergency type~$c$.
We used $\theta_{c} = 1$ for low priority emergencies and $\theta_{c} = 4$ for high priority emergencies.
\item
In \eqref{costalloc}, $M_{ac}$ is the cost (the unit of this cost is the time unit used to measure the response time~$t$ in \eqref{costalloc} and \eqref{penbtct}) of assigning an ambulance of type~$a$ to an emergency of type~$c$, given in Table \ref{tbl:compatibility_matrix}.
\end{itemize}

\begin{table}[H]
\centering
\begin{tabular}{|c|c|c|c|c|}
\hline
   & 1:  High, ALS pref.  & 2: Low, ALS pref. & 3: High, no pref.  & 4: Low, no pref.\\ \hline
ALS & 0 & 0 & 1500 & 1500 \\ \hline
BLS & 6000 & 6000 & 0 & 0 \\ \hline
\end{tabular}
\caption{Quality of care coefficients $M_{ac}$.
The columns correspond to the emergency types while the rows correspond to the ambulance types.}
\label{tbl:compatibility_matrix}
\end{table}

The code was implemented  in C++17 and compiled using GCC~13.1.
The optimization problems were solved using Gurobi version 11.0 and the linear systems were solved using the Eigen3 library.
The computational tests were performed on a computer with AMD Ryzen 5 2600 3.4 Ghz CPU, 16GB of RAM, in a Ubuntu 22.04 OS.

Table~\ref{run_times_markov_heuristic} reports the computation times for the Markov preparedness policy described in Sections~\ref{sec:selection} and~\ref{sec:reassignment}.
The first column shows the number of ambulances which is equal to the number of ambulance stations.
Half the ambulances were ALS ambulances and half were BLS ambulances.
The next four columns report the computation times for solving the optimization problems~\eqref{eqn:linear selection policy objective}--\eqref{eqn:linear selection each emergency at most one dispatch} and \eqref{eqn:linear reassignment policy objective}--\eqref{eqn:linear reassignment policy at most one ambulance} during the simulation (column ``min'' shows the minimum computation time, column ``mean'' shows the mean computation time, column ``q0.9'' shows the $0.9$~quantile of the empirical computation time distribution, and column ``max'' shows the maximum computation time).
The last column shows the total up-front computation time needed to compute the stationary probabilities of the Markov chain and the preparedness metric, using the conjugate gradient method.
Note that most computational effort is needed once, before application of the policy in operations, and that the allocation decisions are computed very fast.
When the number of ambulances exceeds~$16$, then the up-front computation time becomes very large, so we use the following approach.
First, recall that the up-front computations produce the steady state cost rates $\bar{\psi}^{m}_{b}$ for each station $b \in \mathcal{B}$, where $m$ denotes the vector of currently available ambulances at station~$b$.
For example, if there are $M_{1}$~ALS ambulances and $M_{2}$~BLS ambulances in the entire system, then the set of possible values of $m$ is $\mathcal{M}(M_{1},M_{2}) \defi \{0,1,\ldots,M_{1}\} \times \{0,1,\ldots,M_{2}\}$.
For each station $b \in \mathcal{B}$, we compute up-front $\bar{\psi}^{m}_{b}$ for all $m \in \mathcal{M}\left(\min\{8,M_{1}\},\min\{8,M_{2}\}\right)$ only, to keep the up-front computation time reasonable.
If, at any point in time, there are $m_{1}$~ALS ambulances and $m_{2}$~BLS ambulances available at station~$b$, then we use the value $\bar{\psi}^{(\min\{8,m_{1}\},\min\{8,m_{2}\})}_{b}$ as input to the optimization problems~\eqref{eqn:linear selection policy objective}--\eqref{eqn:linear selection each emergency at most one dispatch} and \eqref{eqn:linear reassignment policy objective}--\eqref{eqn:linear reassignment policy at most one ambulance}.
It almost never happens that the number of ALS ambulances or the number of BLS ambulances available at a single station is more than~$8$, so most of the time $(\min\{8,m_{1}\},\min\{8,m_{2}\}) = (m_{1},m_{2})$.
In addition, when $m_{1} > 8$ or $m_{2} > 8$, the incremental value of the additional number of ambulances at station~$b$ beyond~$8$ is very small, so $\bar{\psi}^{(\min\{8,m_{1}\},\min\{8,m_{2}\})}_{b}$ is close to $\bar{\psi}^{(m_{1},m_{2})}_{b}$.

\begin{table}[hbtp]
\centering
\begin{tabular}{@{}cccccc@{}}
\toprule
\multicolumn{1}{c}{\# Ambulances}  & min [ms] & mean [ms] & q0.9 [ms] & max [ms] &  Markov EMS [s] \\ \midrule
6  & 0.02    & 0.74     & 1.31      & 2.10                  & 0.70         \\
8  & 0.019   & 0.43     & 1.11      & 1.79                  & 21.6        \\
10 & 0.02    & 0.19     & 0.66      & 1.52                  & 345      \\ 
12 & 0.06    & 0.27     & 0.41      & 2.19                  & 29100     \\
14 & 0.11    & 0.38     & 0.52      & 2.33                  & 82700      \\
16 & 0.33    & 4.11     & 4.56      & 16.39                 & 346000    \\         
\bottomrule
\end{tabular}
\caption{Computation times for the Markov preparedness policy described in Sections~\ref{sec:selection} and~\ref{sec:reassignment}.}
\label{run_times_markov_heuristic}
\end{table}

We compare the Markov preparedness (MP) policy with the following $9$~policies from the literature:
\begin{enumerate}
\item
The classic closest available ambulance policy (CA) for ambulance selection.
For ambulance reassignment, when there are no emergencies in queue, then the ambulance goes to the closest station.
When there are emergencies in queue, then the ambulance is dispatched to the oldest emergency in queue (the first-come-first-served policy).
\item
For ambulance selection, the method in \citet{ande:07} that dispatches the ambulance~$a$ that maximizes $\min\{\psi^{-a}_{\ell} : \ell \in \mathcal{L}\}$, with $\gamma^{a} = 1$ for all $a \in \mathcal{A}$.
For ambulance reassignment, when there are no emergencies in queue, then the newly available ambulance is sent to a station selected as follows: for any station~$b \in \mathcal{B}$ and zone~$\ell \in \mathcal{L}$, let $\psi^{b+}_{\ell}$ denote the preparedness metric of zone~$\ell$ after the newly available ambulance is added to station~$b$, with $\gamma^{a} = 1$ for all $a \in \mathcal{A}$.
Then the ambulance is sent to the ambulance station~$b$ that maximizes $\min\{\psi^{b+}_{\ell} : \ell \in \mathcal{L}\}$.
If multiple stations attain the maximum, then the ambulance is sent to the maximizing station that is closest to its current position.
When there are emergencies in queue, then the ambulance is dispatched to the oldest emergency in queue.
\item
For ambulance selection, the method in \citet{lees:11} that dispatches the ambulance~$a$ that maximizes $\min\{\psi^{-a}_{\ell} : \ell \in \mathcal{L}\} / t^{a}_{e}$, with $\gamma^{a} = 1$ for all $a \in \mathcal{A}$.
For ambulance reassignment, when there are no emergencies in queue, then the ambulance goes to the closest station.
When there are emergencies in queue, then the ambulance is dispatched to the closest emergency in queue.
\item
The combined policy in \citet{lees:14}, that includes ideas from the method in \citet{lees:12,lees:13}. When there are no emergencies in queue, then the ambulance goes to the closest station.
\item
For ambulance selection, the method in \citet{lees:17} with $w = 1$.
For ambulance reassignment, when there are no emergencies in queue, then the ambulance~$a$ is dispatched to the station $b \in \mathcal{B}$ that maximizes $\min\{\psi^{b+}_{\ell} : \ell \in \mathcal{L}\} / t^{a}_{b}$, where as before $\psi^{b+}_{\ell}$ denotes the preparedness metric of zone~$\ell$ after adding ambulance~$a$ to station~$b$ (also with $\gamma^{a} = 1$ for all $a \in \mathcal{A}$), and $t^{a}_{b}$ denotes the travel time from the current location of ambulance~$a$ to station~$b$.
When there are emergencies in queue, then the method in \citet{lees:12,lees:13} with weighted degree centrality measure $c_{e} = \sum_{e' \in \mathcal{Q} \setminus \{e\}} 1 / \left(1 + t_{e,e'}\right)$ and $w =$ the probability that the ambulance does not transport the patient of emergency~$e$ to a hospital is used.
\item
For ambulance selection, the Heuristic-Cross dispatch rule in \citet{mayo:13}.
For ambulance reassignment, when there are no emergencies in queue, then the ambulance is sent back to its home station.
When there are emergencies in queue, then the ambulance is dispatched to the oldest emergency in queue. 
\item
For ambulance selection, the method in \citet{band:14}.
For ambulance reassignment, when there are no emergencies in queue, then the ambulance is sent back to its home station.
When there are emergencies in queue, then the ambulance is dispatched to the oldest emergency in queue. 
\item
The policy in \citet{jagt:17a}.
Specifically, for ambulance selection, the method based on the MEXCLP preparedness metric is used.
For ambulance reassignment, when there are no emergencies in queue, then the ambulance is sent back to its home station.
When there are emergencies in queue, then the ambulance is dispatched to the oldest emergency in queue.
\item
The policy proposed in \citet{carv:25}, for the setting in which each emergency only needs one ambulance. 
\end{enumerate}

We compare results for two spatial discretizations of the city: a discretization into $10 \times 10$ rectangular zones and a hexagonal discretization.
The time of the week was partitioned into $30$~minute time intervals.
A nonhomogeneous Poisson arrival process for each emergency type was calibrated for both spatial discretizations using historical data and maximizing a regularized likelihood function, as described in \citet{laspatedpaper,laspatedmanual}.
Results were produced for two distance metrics and associated travel times: (a)~great circle distances and a constant speed of 60km/h for the ambulances and (b)~travel distances along the streets of Rio de Janeiro, with travel times that vary depending on the type of road being traversed, as provided by the C++ OpenStreetMap library.
We solved the system of equations \eqref{eqn:balance equations}--\eqref{eqn:probability scaling} for up to 16 ambulances, and computed the resulting Markov preparedness ambulance selection and ambulance reassignment policies described in Sections~\ref{sec:selection} and~\ref{sec:reassignment} respectively.

\ignore{
The heuristics were applied both in rollout mode and in nonrollout mode (see \citealt{guiguesetal2026} for details on how the heuristics are used in rollout mode, solving two stage stochastic problems each time a decision has to be made).
Using OpenStreetMap significantly increased the run time of the simulations, so we did not evaluate the heuristics in rollout mode in the experiments that used OpenStreetMap to compute travel times of ambulances along the streets of the city.}

\begin{figure}
    \centering
    \includegraphics[width=0.9\linewidth]{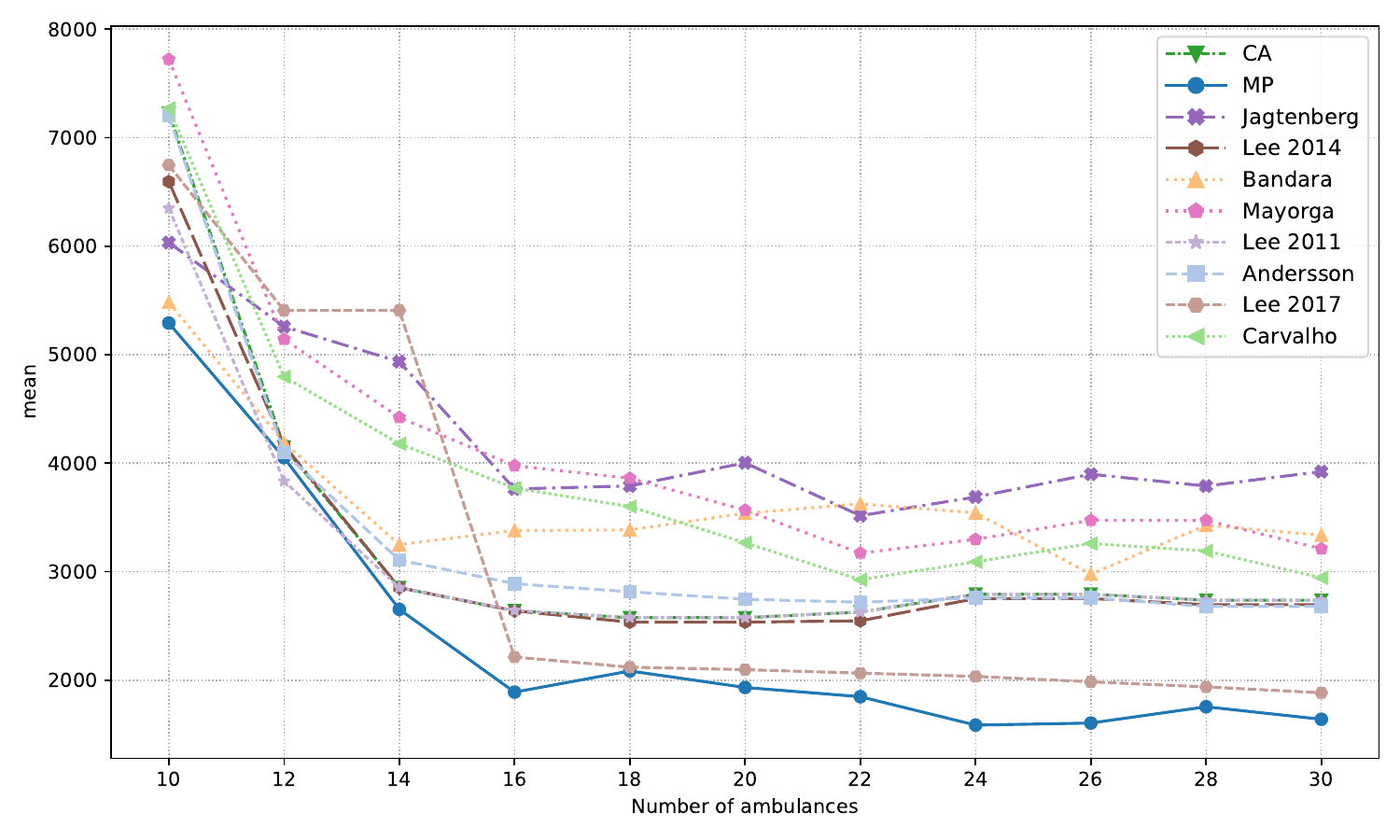}
    \caption{Simulated mean allocation costs for the policies, using great circle distances and rectangular discretization.}
    \label{fig:rectangular_allocation_costs}
\end{figure}

\begin{figure}
    \centering
    \includegraphics[width=0.9\linewidth]{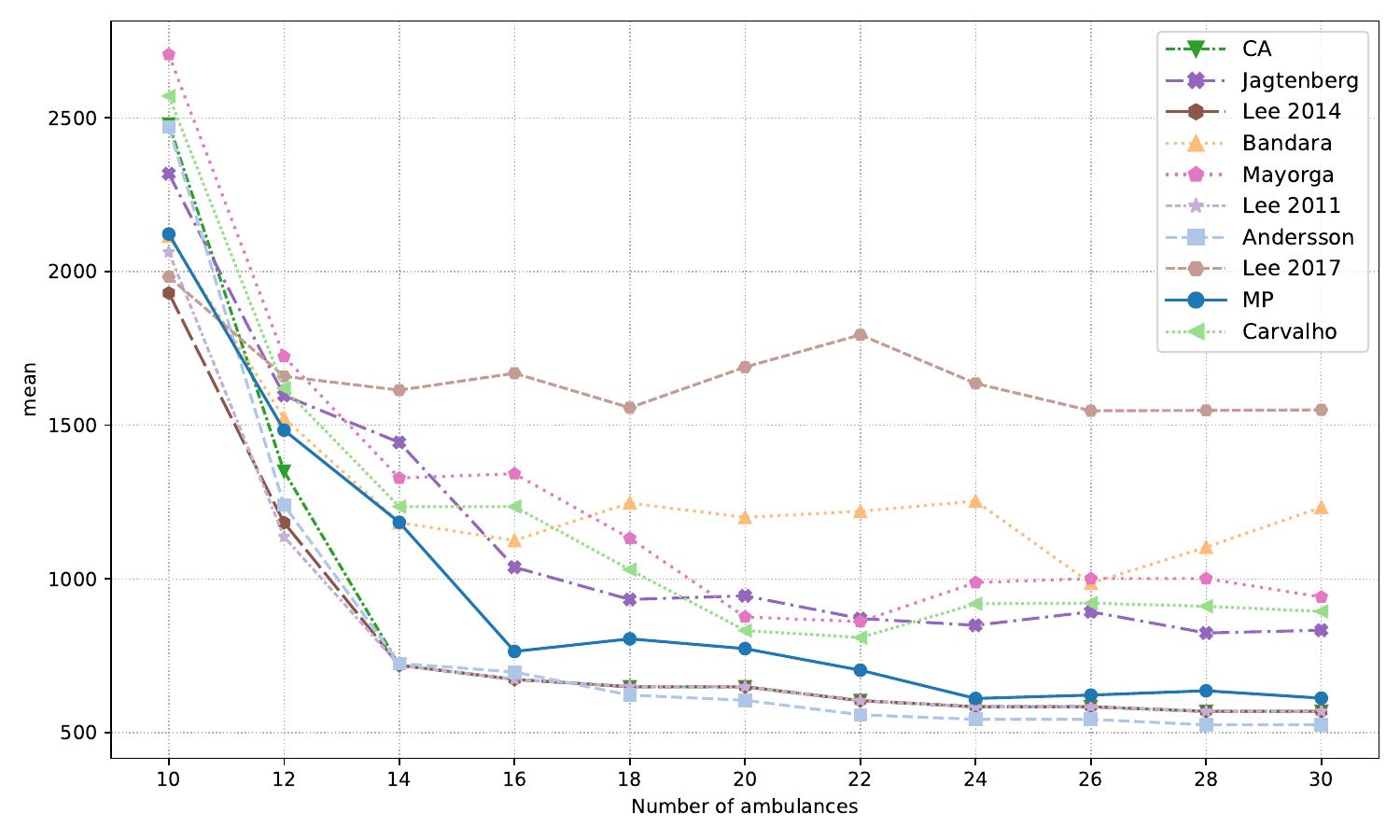}
    \caption{Simulated mean response times for the policies, using great circle distances and rectangular discretization.}
    \label{fig:placeholder3}
\end{figure}

\begin{figure}
    \centering
    \includegraphics[width=0.9\linewidth]{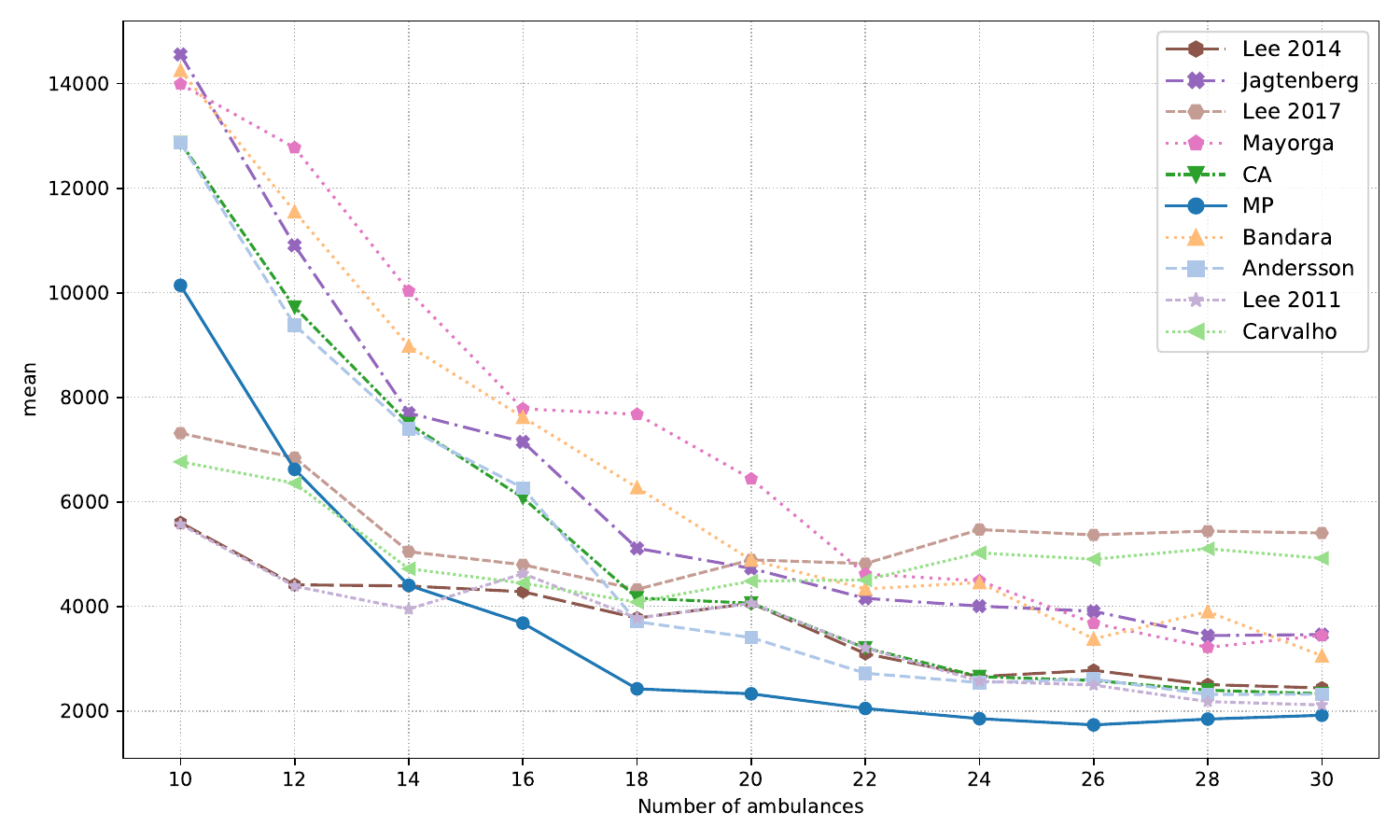}
    \caption{Simulated mean allocation costs for the policies, using great circle distances and hexagonal discretization.}
    \label{fig:hex7_allocation_costs}
\end{figure}

\begin{figure}
    \centering
    \includegraphics[width=0.9\linewidth]{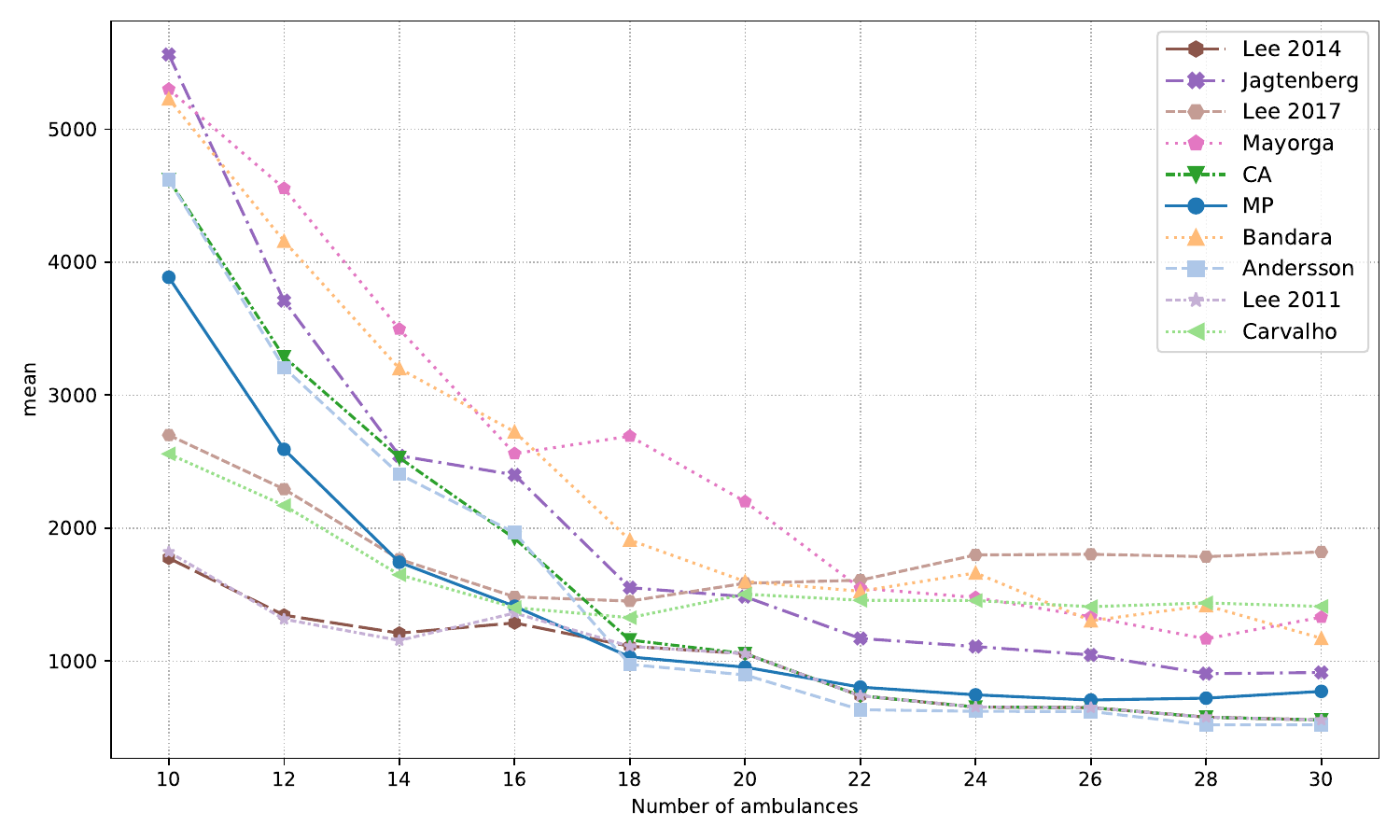}
    \caption{Simulated mean response times for the policies, using great circle distances and hexagonal discretization.}
    \label{fig:hex7_response_times}
\end{figure}

\begin{figure}
    \centering
    \includegraphics[width=0.9\linewidth]{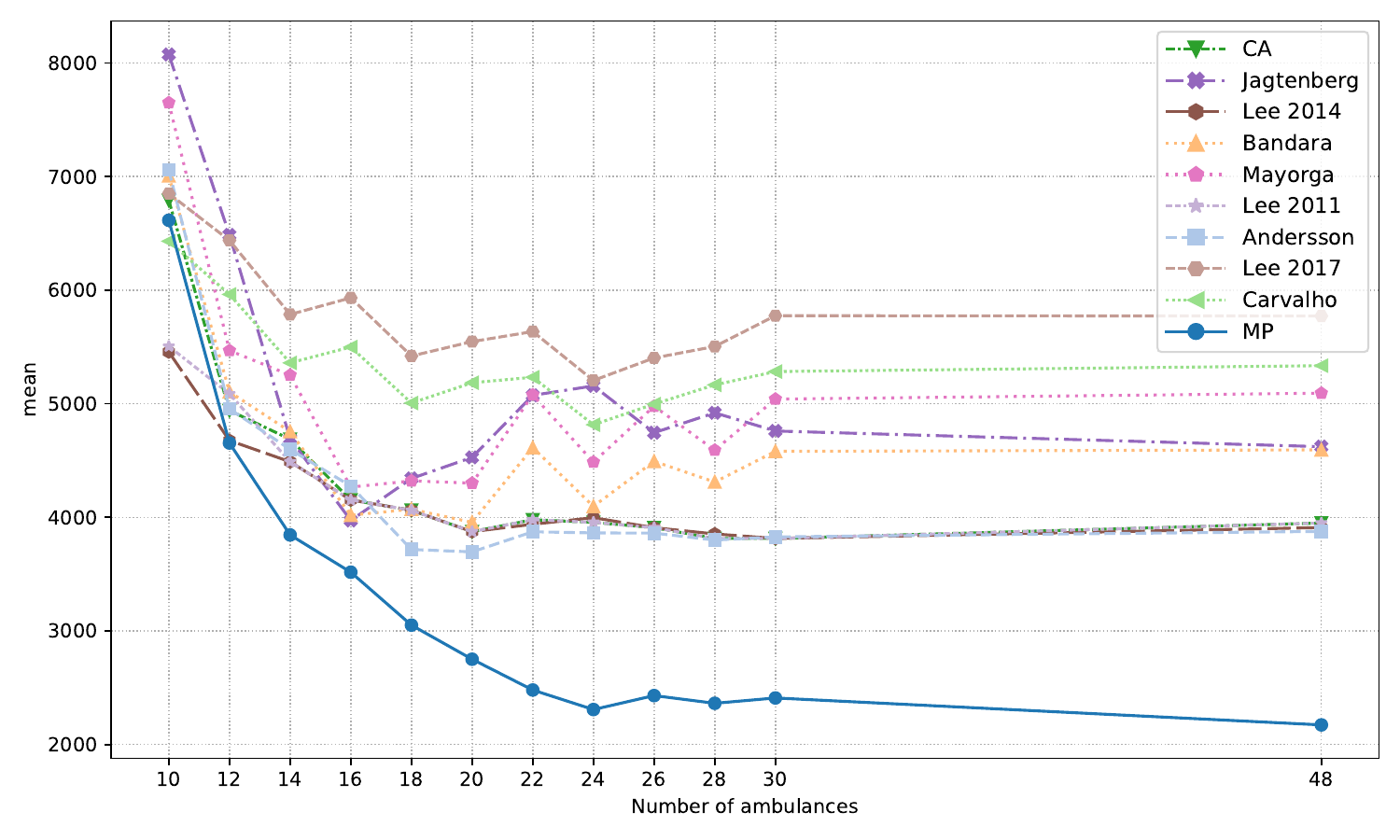}
    \caption{Simulated mean allocation costs for the policies, using rectangular discretization and travel times along the street network of Rio de Janeiro.}
    \label{fig:street_rect_allocation_costs}
\end{figure}

\begin{figure}
    \centering
    \includegraphics[width=0.9\linewidth]{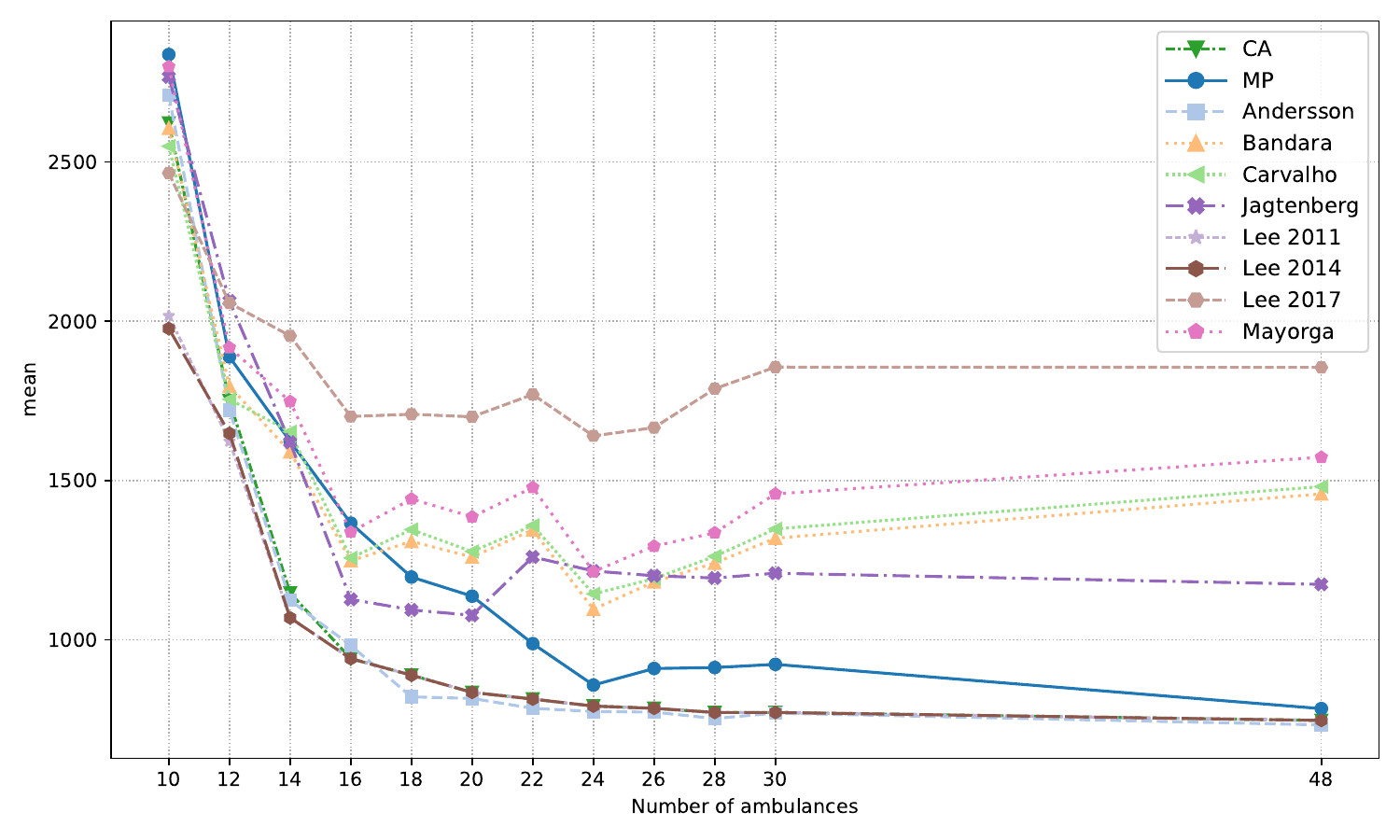}
    \caption{Simulated mean response times for the policies, using rectangular discretization and travel times along the street network of Rio de Janeiro.}
    \label{fig:street_rect_response times}
\end{figure}

\begin{figure}
    \centering
    \includegraphics[width=0.9\linewidth]{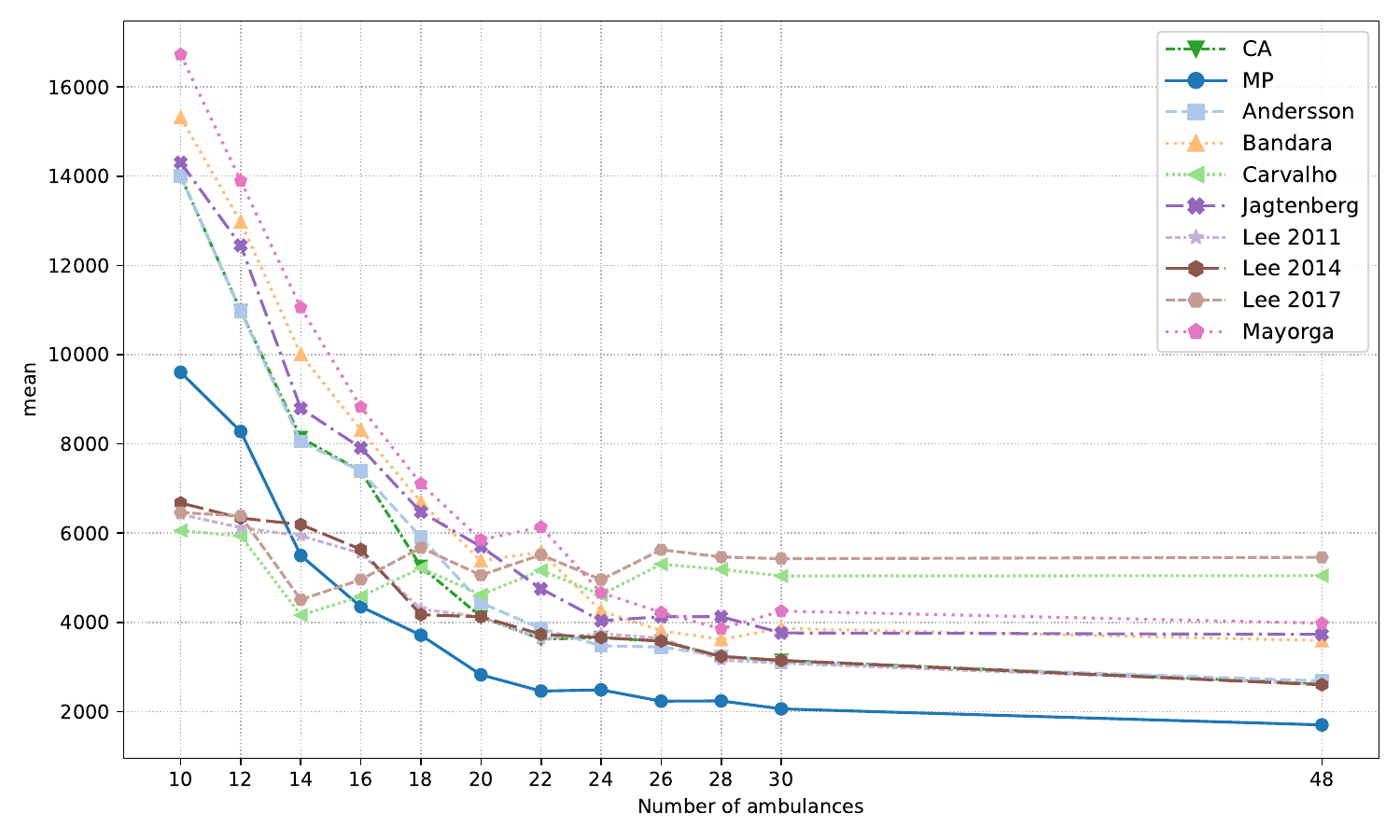}
    \caption{Simulated mean allocation costs for the policies, using hexagonal discretization and travel times along the street network of Rio de Janeiro.}
    \label{fig:street_hex_allocation_costs}
\end{figure}

\begin{figure}
    \centering
    \includegraphics[width=0.9\linewidth]{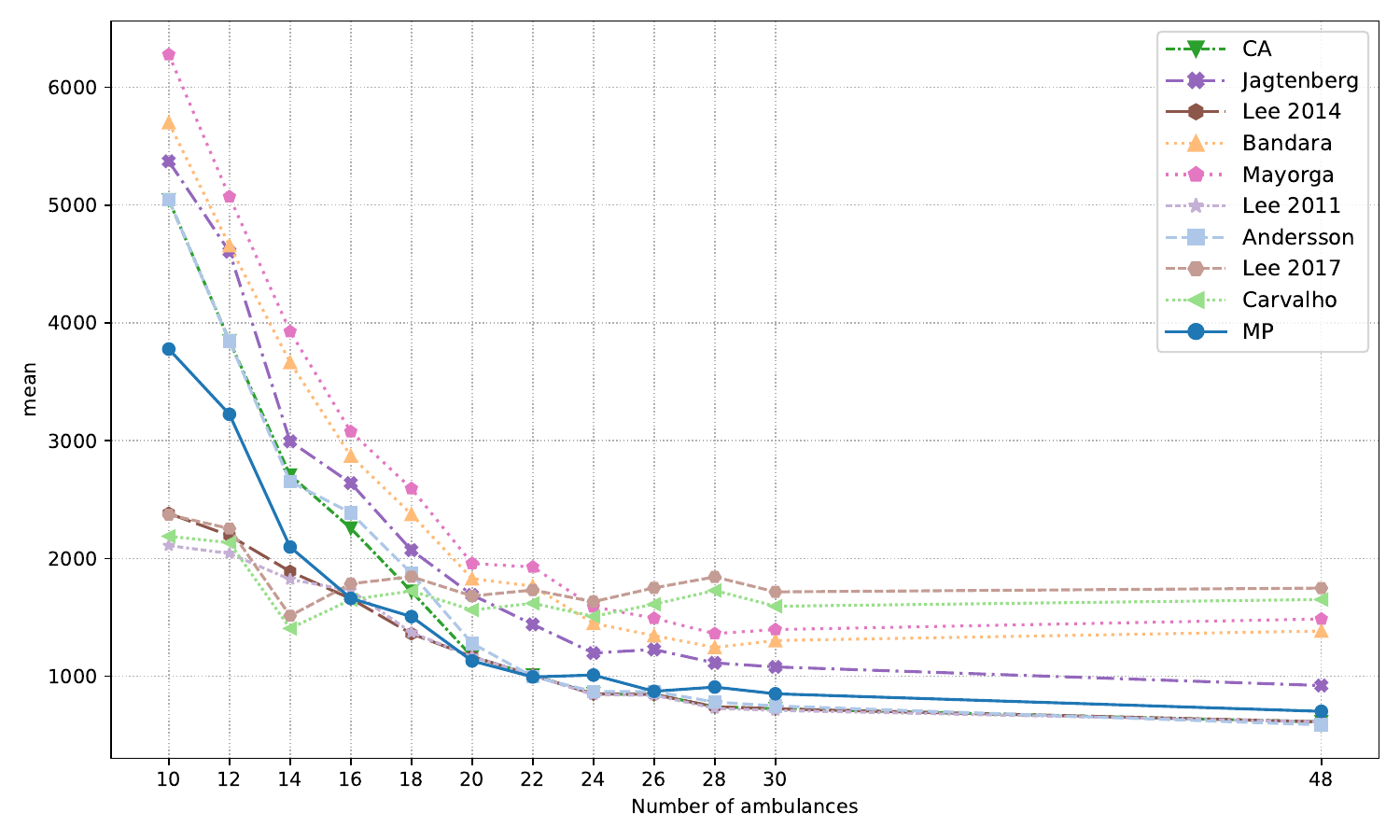}
    \caption{Simulated mean response times for the policies, using hexagonal discretization and travel times along the street network of Rio de Janeiro.}
    \label{fig:street_hex_response times}
\end{figure}

\begin{figure}
    \centering
    \includegraphics[width=0.9\linewidth]{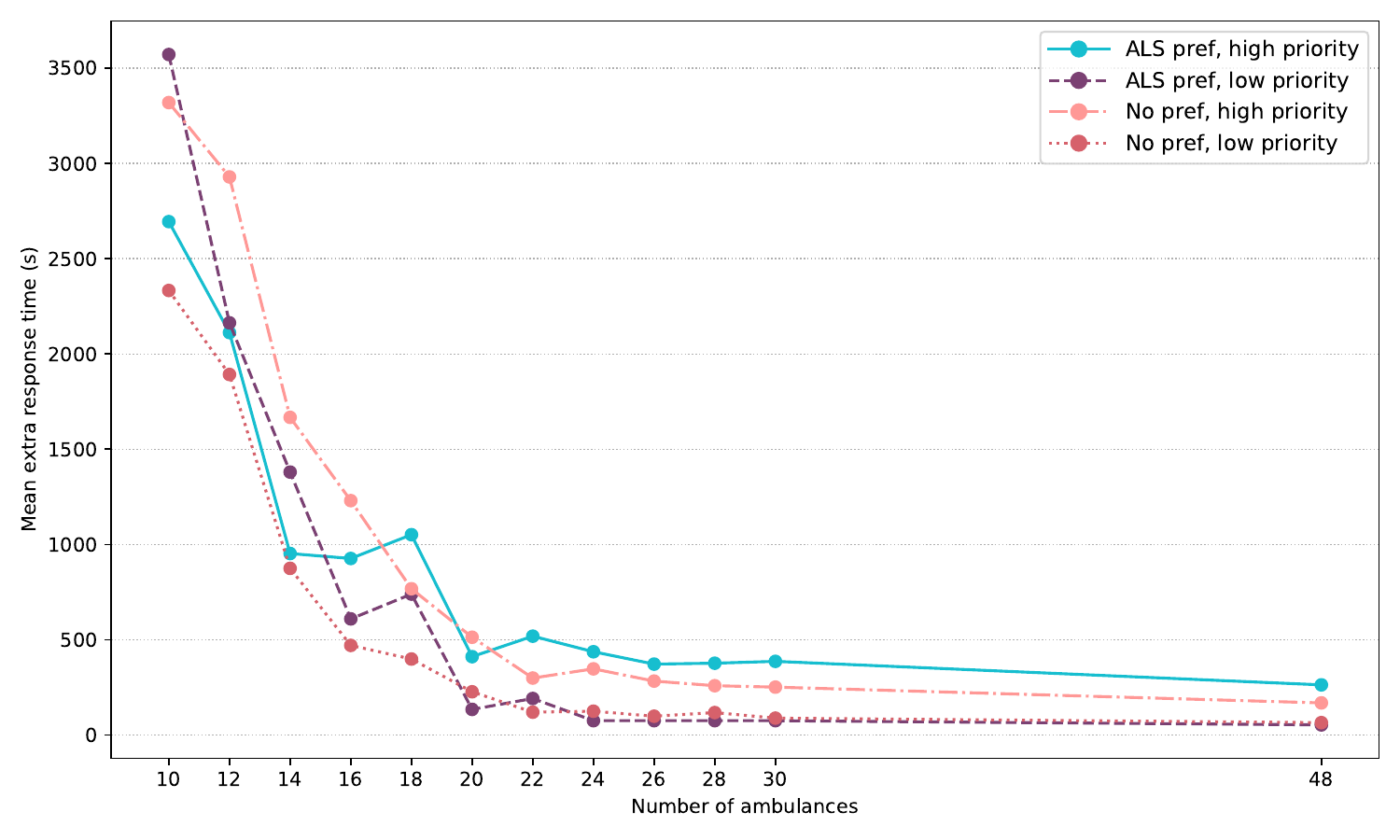}
    \caption{Simulated mean extra response times for each emergency type, relative to a target of 600 seconds for high-priority emergencies and 1200 seconds for low-priority emergencies.}
    \label{fig:mean_extra_response_time}
\end{figure}

The results are reported in Figures~\ref{fig:rectangular_allocation_costs}--\ref{fig:mean_extra_response_time}.
We show both the response times and the allocation costs as given by function \eqref{costalloc} with time penalization function given by \eqref{penbtct}.
Specifically, Figure~\ref{fig:rectangular_allocation_costs} plots the simulated mean allocation costs for the policies using great circle distances and rectangular space discretization.
The MP heuristic provides the best mean allocation costs for all values of the number of ambulances, except for $12$~ambulances in which case the heuristic of \citep{ande:07} is slightly better.
The allocation costs naturally tend to decrease when the number of ambulances increases and flattens off after a sufficient number of ambulances is provided to the system.
Figure~\ref{fig:placeholder3} shows the simulated mean response times using great circle distances and rectangular discretization.
The mean response times of the MP policy are among the best mean response times, even though the MP policy is based on optimization problems~\eqref{eqn:linear selection policy objective}--\eqref{eqn:linear selection each emergency at most one dispatch} and \eqref{eqn:linear reassignment policy objective}--\eqref{eqn:linear reassignment policy at most one ambulance} in which objective coefficients $r(a,i)$ represent allocation costs and not just response times.
Figure~\ref{fig:hex7_allocation_costs} plots the simulated mean allocation costs, and Figure~\ref{fig:hex7_response_times} plots the simulated mean response times, both with great circle distances and hexagonal discretization.
The conclusions are similar to the conclusions described above with rectangular discretization.

Figures~\ref{fig:street_rect_allocation_costs} and~\ref{fig:street_rect_response times} show the simulated mean allocation costs and simulated mean response times respectively, using rectangular spatial discretization and ambulances traveling along the streets of the city.
Figures~\ref{fig:street_hex_allocation_costs} and~\ref{fig:street_hex_response times} show the same performance metrics, using hexagonal discretization.
In most cases, the MP policy is the best in terms of allocation costs and is among the best in terms of response times, except when there is a small number of ambulances.
The relatively poorer performance of the MP policy when there is a small number of ambulances can be explained as follows.
Recall that the Markov model produces the steady state cost rates $\bar{\psi}^{m}_{b}$ for each station~$b$ that are used as input for the optimization problems~\eqref{eqn:linear selection policy objective}--\eqref{eqn:linear selection each emergency at most one dispatch} and \eqref{eqn:linear reassignment policy objective}--\eqref{eqn:linear reassignment policy at most one ambulance} that produce the decisions for the MP policy.
The Markov model assumes that stations operate independently, whereas the optimization problems consider the entire system, and thus the optimization problems can send an ambulance from one station to an emergency at a location usually served from another station, or send an available ambulance close to one station to another station.
Thus, the Markov chain is a simplified model of what happens when the MP policy is used to manage the ambulances.
In a well-supplied and well-operated EMS, most of the time ambulances should travel from a nearby station to serve an emergency, in which case we expect the Markov chain to be a reasonably accurate model of what happens under the MP policy.
Things may go wrong if the EMS is not well-supplied or not well-operated.
The MP policy strives to make the EMS well-operated, so the main case when the Markov model may not be accurate is when an EMS is not well-supplied, that is, ambulance supply is small relative to demand so that ambulances often have to be dispatched from far away to serve an emergency.
Thus, when the number of ambulances is small, the preparedness metric produced by the Markov model may not be accurate, and the MP policy may not perform well.
In such a setting the main problem is not the policy used for EMS operations, but rather lack of EMS resources.

Figure~\ref{fig:mean_extra_response_time} shows the simulated mean extra response time under the MP policy as a function of the number of ambulances, for each emergency type, using street travel times and hexagonal discretization.
Here, the extra response time of an emergency is the response time in excess of $10$~minutes for high-priority emergencies and $20$~minutes for low-priority emergencies.

In summary, the main conclusion from the numerical results is that the MP policy performs better than the other policies in terms of the user's chosen performance measure (recall that the user can choose any cost function $\phi$ and that the preparedness metric is based on the steady-state average of $\phi$), except when the supply of ambulances is very small relative to the demand.
In addition, the MP policy performs among the best in terms of expected response time, except when the supply of ambulances is small relative to the demand.
The situation when the supply of ambulances is too small for the demand is difficult, and there does not seem to be a single policy that consistently performs well relative to the others when the EMS is under-supplied, with the policies of \citet{lees:11,lees:14} appearing to be the safest choices in such a situation.

\section{Source Code and Data}
\label{sec:data}

The source code and data for all the policies considered in this work (the Markov preparedness policy and the policies from the literature) are available in the "esma2" directory on Github at  {\url{https://github.com/vguigues/Heuristics_Dynamic_Ambulance_Management}} and on Zenodo at {\url{DOI: 10.5281/zenodo.18487224}}.
The code requires the following libraries:

\begin{itemize}
    \item Boost, available at {\url{https://boost.org}},
    \item xtl and xtensor, both available at \url{https://github.com/xtensor-stack},
    \item fmt, available at {\url{https://github.com/fmtlib}},
    \item Gurobi, available at {\url{https://gurobi.com}}
    \item OSRM library, available at {\url{https://github.com/Project-OSRM/osrm-backend}}
\end{itemize}

The code was tested on a Linux Ubuntu operating system.
The Boost library can be acquired using:

\begin{verbatim}
    sudo apt install libboost-all-dev
\end{verbatim}

The libraries available at GitHub can be installed via CMake with the following commands, from the root directory of each library source code:

\begin{verbatim}
    mkdir build
    cd build
    cmake ..
    make
    sudo make install
\end{verbatim}

To compile the simulation code, just run the following commands from the esma2 directory:

\begin{verbatim}
    mkdir build
    cd build
    cmake -DESMA_ENABLE_TRAJECTORIES=OFF 
    -DCMAKE_BUILD_TYPE=Release -DESMA_STREET_TRAVEL=ON ..
    make
\end{verbatim}

The command above will compile the code and use OSRM to compute travel times/distances. If you want to use geodesic distances, just omit the ESMA\_STREET\_TRAVEL option. The ESMA\_ENABLE\_TRAJECTORIES option enables the simulator to write the ambulance trajectories in the trajectories' directory. In this example, we turn off trajectories to improve the simulation performance. 

The Gurobi library should be extracted at the /opt directory and an environment variable \$GRB\_LICENSE must be set  with the path to the Gurobi License. 
After installing the dependencies, run the make command from the source code root directory.
After successfully compiling the source code, the executable can be run as in the following example:

\begin{verbatim}
    ./esma -f test.cfg --amb_setup=rj --n_scenarios=25 --nb_ambulances=6 
    --nb_bases=6
\end{verbatim}

The command above runs 25 scenarios/replications for each policy, both with and without rollout.
(In rollout mode, each time a decision has to be made a two-stage stochastic program is solved in which the second stage objective value is given by the performance under the considered policy.
See \citealt{guiguesetal2026} for more details on how the policies are used in rollout mode.)

The simulation in this example is based on the Rio de Janeiro EMS System rules (parameter amb\_setup) with 6 ambulances (parameter n\_ambulances) and 6 ambulance stations (parameter n\_bases, with a maximum of 34).
If parameter amb\_setup is set to ``us'' then American EMS system rules are used.

The results for each policy are saved in the results directory, for the policies without rollout and with rollout, respectively.
The response times and allocation costs of each policy are saved in files with the following naming convention:
$$
\mbox{setup\_policy\_n\_ambulances\_n\_scenarios.txt},
$$
denoting respectively the EMS rules, the policy, the number of ambulances and the number of scenarios.
The possible names for ``policy'' in the file names are ``dummy\_queue'' for the closest available ambulance policy, ``markov\_preparedness'' for the Markov preparedness policy using the models of Sections~\ref{sec:selection} and~\ref{sec:reassignment}, ``preparedness'' for the policy of \citet{ande:07}, ``prep2'' for the policy of \citet{lees:11}, ``district'' for the policy of \citet{mayo:13}, ``ordered'' for the policy of \citep{band:14}, ``coverage'' for the policy of \citep{jagt:17a}, ``centrality'' for the policy of \citet{lees:14}, ``dist\_centrality'' for the policy of \citet{lees:17}, and ``tipat'' for the policy of \citet{carv:25}.
Each result file contains, for each scenario, the number of emergencies \(N\) in the scenario, followed by \(N\) lines containing: the index of the ambulance that served the emergency, the response time for the emergency, its allocation cost, and the instant the ambulance finished service for that emergency.

\section{Conclusion}

In this paper, we introduced a new policy for the operation of an ambulance fleet under uncertainty based on a new preparedness metric.
This policy was compared with $9$~policies proposed in the literature using a simulation based on data of the Rio de Janeiro EMS.
In most cases, this policy performed better than the $9$~policies in terms of mean allocation costs, and this policy also performed among the best of the $9$~policies in terms of mean response time.

As future work, we intend to adapt the policy to consider a more diverse fleet of vehicles (including automobiles, ALS and BLS units, motorcycles, drones, and helicopters), taking into account the differences in capacity of the vehicles in terms of the number of patients and equipment that can be transported.

\section*{Data Availability Statement}

The source code and data used to produce the results of this paper are available on GitHub and Zenodo.
See Section~\ref{sec:data} for instructions on how to download and use the source code and data.

\section*{Declaration of interests} 

There are no competing interests.

\section*{Declaration of funding} 

No funding was received.

\addcontentsline{toc}{section}{References}
\printbibliography

\end{document}